\documentclass{article}
\usepackage{amsthm,amsfonts,amsmath,amscd,amssymb,latexsym,epsfig}
\title{Wall crossing for symplectic vortices\\
       and quantum cohomology}

\author{Kai~Cieliebak\thanks{Supported by National Science
Foundation Grant DMS--0072267}\\
Ludwig--Maximilans--Universit\"at M\"unchen
\and
Dietmar~A.~Salamon\\
ETH--Z\"urich}

\date{18 September 2002}

\newtheorem{theorem}{Theorem}[section]
\newtheorem{corollary}[theorem]{Corollary}

\newtheorem{lemma}[theorem]{Lemma}
\newtheorem{proposition}[theorem]{Proposition}

\newtheorem{remark}[theorem]{Remark}

\newtheorem{example}[theorem]{Example}
\newcommand{\one}{{{\mathchoice {\rm 1\mskip-4mu l} {\rm 1\mskip-4mu l}
{\rm 1\mskip-4.5mu l} {\rm 1\mskip-5mu l}}}}
\newcommand{\dslash}{/\mskip-6mu/}

\newcommand{\B}{{\mathbb{B}}}
\newcommand{\C}{{\mathbb{C}}}
\newcommand{\E}{{\mathbb{E}}}

\newcommand{\HH}{{\mathbb{H}}}
\newcommand{\LL}{{\mathbb{L}}}
\newcommand{\M}{{\mathbb{M}}}

\renewcommand{\P}{{\mathbb{P}}}

\newcommand{\R}{{\mathbb{R}}}
\renewcommand{\SS}{{\mathbb{S}}}
\newcommand{\T}{{\mathbb{T}}}
\newcommand{\Z}{{\mathbb{Z}}}
\newcommand{\Aa}{{\mathcal{A}}} 
\newcommand{\Bb}{{\mathcal{B}}}
\newcommand{\Dd}{{\mathcal{D}}}
\newcommand{\Ee}{{\mathcal{E}}}
\newcommand{\Ff}{{\mathcal{F}}}
\newcommand{\Gg}{{\mathcal{G}}} 

\newcommand{\Ii}{{\mathcal{I}}}
\newcommand{\Jj}{{\mathcal{J}}}

\newcommand{\Ll}{{\mathcal{L}}} 
\newcommand{\Mm}{{\mathcal{M}}} 

\newcommand{\Pp}{{\mathcal{P}}}

\newcommand{\Rr}{{\mathcal{R}}}
\newcommand{\Ss}{{\mathcal{S}}}

\newcommand{\Ww}{{\mathcal{W}}}
\newcommand{\Xx}{{\mathcal{X}}}
\newcommand{\Yy}{{\mathcal{Y}}}

\newcommand{\coker}{{\rm coker }} 
\newcommand{\im}{{\rm im }} 
\newcommand{\SPAN}{{\rm span }} 
\newcommand{\id}{{\rm id}} 

\newcommand{\codim}{{\rm codim}} 
\newcommand{\INDEX}{{\rm index}} 
\newcommand{\IND}{{\mathcal{IND}}} 

\newcommand{\IM}{{\rm Im }} 
\newcommand{\RE}{{\rm Re}} 
\newcommand{\Lie}{{\rm Lie}} 
\newcommand{\Aut}{{\rm Aut}} 
\newcommand{\Vect}{{\rm Vect}} 
\newcommand{\Vol}{{\rm Vol}} 
\newcommand{\Res}{{\rm Res}} 
\newcommand{\End}{{\rm End}} 
\newcommand{\GW}{{\rm GW}} 
\newcommand{\QH}{{\rm QH}}

\newcommand{\W}{{\rm W}}
\newcommand{\w}{{\rm w}}

\newcommand{\point}{{\rm pt}}

\newcommand{\e}{{\rm e}}

\newcommand{\eff}{{\rm eff}}

\newcommand{\dvol}{{\rm dvol}}

\newcommand{\al}{{\alpha}}

\newcommand{\eps}{{\varepsilon}}

\renewcommand{\i}{{\iota}}
\newcommand{\la}{{\lambda}}

\newcommand{\om}{{\omega}}

\newcommand{\Om}{{\Omega}}
\newcommand{\G}{{\rm G}}
\newcommand{\HG}{{\rm H}}

\newcommand{\ET}{{\rm ET}}
\newcommand{\BT}{{\rm BT}}

\newcommand{\U}{{\rm U}}

\renewcommand{\u}{{\mathfrak u}}

\newcommand{\g}{{\mathfrak g}} 
\renewcommand{\tt}{{\mathfrak t}} 
\newcommand{\Cinf}{C^{\infty}}

\newcommand{\PD}{{\rm PD}}

\newcommand{\Pic}{{\rm Pic}}

\newcommand{\inner}[2]{\langle #1, #2\rangle}

\def\NABLA#1{{\mathop{\nabla\kern-.5ex\lower1ex\hbox{$#1$}}}}
\def\Nabla#1{\nabla\kern-.5ex{}_{#1}}
\def\Tabla#1{\Tilde\nabla\kern-.5ex{}_{#1}}
\renewcommand{\Tilde}{\widetilde}

\newcommand{\p}{{\partial}}

\newcommand{\IMP}{\Longrightarrow}

\newcommand{\INTO}{\hookrightarrow}
\newcommand{\TO}{\longrightarrow}
\renewcommand{\la}{\langle}
\newcommand{\ra}{\rangle}
\newcommand{\ch}{{\rm ch}}
\newcommand{\td}{{\rm td}}

\begin{document}

\maketitle


\begin{abstract}
We derive a wall crossing formula for the symplectic 
vortex invariants of toric manifolds. As an application, 
we give a proof of Batyrev's formula for the quantum
cohomology of a monotone toric manifold with minimal 
Chern number at least two. 
\end{abstract}



\section{Introduction}

Let $T$ be a torus of dimension $k$,
denote by $\tt$ its Lie algebra, by
$$
\Lambda := \left\{\xi\in\tt\,|\,\exp(\xi)=1\right\}
$$
the integer lattice, and by 
$$
\Lambda^* := \left\{\w\in\tt^*\,|\,\inner{\w}{\xi}\in\Z\text{ for
}\xi\in\Lambda\right\}
$$
the dual lattice. Suppose $T$ acts diagonally on $\C^n$.
The action is determined by $n$ homomorphisms $\rho_\nu:T\to S^1$,
$\nu=1,\dots,n$. We write each homomorphism $\rho_\nu$ in the form
$$
\rho_\nu(\exp(\xi)) = e^{-2\pi i\inner{\w_\nu}{\xi}},\qquad
\w_\nu\in\Lambda^*.
$$
The moment map of this action, with respect to the standard symplectic 
form on $\C^n$, is given by
\begin{equation}\label{eq:mu}
\mu(x) = \pi\sum_{\nu=1}^n\left|x_\nu\right|^2\w_\nu
\end{equation}
for $x=(x_1,\dots,x_n)\in\C^n$. We assume throughout
that $\mu$ is proper and that the vectors $\w_\nu$ span
the space $\tt^*$.  In~\cite{CGMS} we have defined 
invariants
\begin{equation}\label{eq:m}
     \Phi^{\rho,\tau}_\lambda:S^m(\tt^*)\to\R,\qquad
     m := n - \dim\,T + \sum_{\nu=1}^nd_\nu\ge 0,
\end{equation}
by counting solutions of the genus zero symplectic vortex equations 
(see Section~\ref{sec:vortex}).  Here 
$\lambda\in\Lambda$, 
$
     d_\nu:=\inner{\w_\nu}{\lambda},
$
$\tau$ is a regular value of the moment map, 
and $S^m(\tt^*)$ denotes the space of real valued 
polynomials of degree $m$ on $\tt$.  Note that $S^*(\tt^*)$
is canonically isomorphic to the cohomology $H^*(BT;\R)$ of the
classifying space $BT=ET/T$.
The isomorphism takes $\w_\nu\in\tt^*$
to the first Chern class of the bundle 
$ET\times_{\rho_\nu}\C\to BT$.  The invariant  
$\Phi^{\rho,\tau}_\lambda$ takes rational values on 
integral cohomology classes.  These correspond to
polynomials that map the lattice to the integers. 

An element $\tau\in\tt^*$ is a singular value of $\mu$
if and only if it can be expressed as a positive linear 
combination of at most $k-1$ of the vectors $\w_\nu$. 
The set of singular values is a disjoint union of open cones of 
codimensions $1$ to $k$. A cone of codimension $j$ is called a 
{\bf wall of codimension $j$}.  

Let $\tau_0\in\tt^*$ be an element of a wall of codimension one,
$\tau_1\in\tt^*$ be transverse to the wall at $\tau_0$, 
and $e_1\in\Lambda$ be the unique primitive lattice vector 
that is orthogonal to the wall at $\tau_0$
and satisfies $\inner{\tau_1}{e_1}>0$. Denote by $T_1\subset T$ the 
subtorus generated by $e_1$ and by $\tt_1$ its Lie algebra. 
Let 
$$
      I:=\{\nu\;|\;\inner{\w_\nu}{e_1}=0\}.
$$ 
The action $\rho$ induces an action $\rho_0$ of the quotient torus 
$$
      T_0:=T/T_1
$$ 
on the space 
$$
      \C^I:=\{x\in\C^n\;|\;x_\nu=0\mbox{ for }\nu\notin I\}.
$$
The moment map of this action is the restriction
$\mu_0:=\mu|_{\C^I}:\C^I\to\tt_0^*:=\tt_1^\perp$.
The following {\it wall crossing formula} expresses the difference 
of the invariants on the two sides of the wall as the invariant of 
the reduced problem at $\tau_0$.

\begin{theorem}[Genus Zero Wall Crossing]\label{thm:wall0}
Let $\alpha\in S^*(\tt^*)$, $\lambda\in\Lambda$, and 
$d_\nu:=\inner{\w_\nu}{\lambda}$.  Then for every
sufficiently small positive number $\eps$ we have
$$
      \Phi^{\rho,\tau_0+\eps\tau_1}_\lambda(\al)
      - \Phi^{\rho,\tau_0-\eps\tau_1}_\lambda(\al)
      = \Phi^{\rho_0,\tau_0}_{\lambda_0}(\al_0),
$$
where $\lambda_0$ is the projection of $\lambda$ to 
$\tt_0:=\tt/\tt_1$ and
$$
     \alpha_0(\xi) = \frac{1}{2\pi i}
     \oint\frac{\alpha(\xi+ze_1)}{\prod_{\nu\notin I}
     \inner{\w_\nu}{\xi+ze_1}^{d_\nu+1}}dz.
$$
Here for each $\xi$ the integral is understood over a circle in the 
complex plane enclosing all the poles of the integrand.
\end{theorem}

There is an analogous wall crossing formula for higher genus 
which is formulated in Theorem~\ref{thm:WCN} below.  

Theorem~\ref{thm:wall0} gives rise to an explicit formula for 
the genus zero invariants.  To formulate the result we 
introduce the following notation.  For a tuple of nonnegative
integers $\ell=(\ell_1,\dots,\ell_n)$ denote 
$$
     \w^\ell := \w_1^{\ell_1}\cdots\w_n^{\ell_n}
     \in S^{|\ell|}(\tt^*),\qquad
     |\ell| := \ell_1+\dots+\ell_n.
$$
Given such a tuple $\ell$ and a lattice vector $\lambda$ we introduce 
the set $\Ii_\lambda(\ell)$ of partitions
$
      \{1,\dots,n\} = I_1\cup\cdots\cup I_k
$
that satisfy the following two conditions. 

\begin{description}
\item[(Dimension)]
For every $j\in\{1,\dots,k\}$, the subspace 
$$
      E_j:=\SPAN\{\w_\nu\,|\,\nu\in I_1\cup\cdots\cup I_j\}
      \subset\tt^*
$$
has dimension $j$ and $\w_\nu\notin E_j$ for every 
$\nu\in I_{j+1}\cup\cdots\cup I_k$. 
\item[(Degree)]
For every $j\in\{1,\dots,k\}$,
$$
     \sum_{\nu\in I_j}(\ell_\nu-d_\nu-1)=-1.\qquad
     d_\nu := \inner{\w_\nu}{\lambda}.
$$
\end{description}

\begin{theorem}[Genus Zero Invariants]\label{thm:invariants}
Let $\lambda\in\Lambda$, $d_\nu:=\inner{\w_\nu}{\lambda}$, and 
$\ell$ be an $n$-tuple of nonnegative integers. 

\smallskip
\noindent{\bf (i)}
If $\Ii_\lambda(\ell)=\emptyset$ then 
$\Phi^{\rho,\tau}_\lambda(\w^\ell)=0$ for all $\tau$.

\smallskip
\noindent{\bf (ii)}
Let $J\subset\{1,\dots,n\}$ be a subset with $k$ elements such 
that $\{\w_\nu\,|\,\nu\in J\}$ is a basis of $\tt^*$ and assume
$$
      \ell_\nu=\begin{cases} 
      d_\nu &\mbox{if }\nu\in J, \\ 
      d_\nu+1 &\mbox{if }\nu\notin J.
      \end{cases}
$$ 
If $\tau$ belongs to the cone $C(J)$ spanned by 
$\{\w_\nu\,|\,\nu\in J\}$ then
$$
     \Phi^{\rho,\tau}_\lambda(\w^\ell) 
     = \frac{1}
       {|\det(\inner{\w_\nu}{e_j}_{\nu\in J,j=1,\dots,k})|},
$$
where $e_1,\dots,e_k$ is any basis 
of the lattice $\Lambda$. Otherwise 
$\Phi^{\rho,\tau}_\lambda(\w^\ell)=0$.

\smallskip
\noindent{\bf (iii)}
Let $\lambda'\in\Lambda$ and define $d'_\nu:=\inner{\w_\nu}{\lambda'}$.
If $\ell_\nu+d'_\nu\ge 0$ for every $\nu$ then 
$$
     \Phi^{\rho,\tau}_{\lambda}(\w^\ell)
     = \Phi^{\rho,\tau}_{\lambda+\lambda'}(\w^{\ell+d'}).
$$

\smallskip
\noindent{\bf (iv)}
Assume $d_\nu\ge-1$ for every $\nu$.
Then every element of $S^*(\tt^*)$ is a linear combination of 
monomials $\w^\ell$ that satisfy either~(i) or~(ii).

\smallskip
\noindent{\bf (v)}
Let $J_\ell:=\left\{\nu\,|\,\ell_\nu\le d_\nu\right\}$.
If $\tau\notin C(J_\ell)$ then 
$\Phi^{\rho,\tau}_\lambda(\w^\ell)=0$.
\end{theorem}

\noindent{\bf Remark.}
Assertions~(i), (ii) and~(iv) can be used to compute 
the genus zero invariants whenever $d_\nu\ge -1$. 
This restriction can be removed by using~(iii). 

\medskip
Now assume that $T$ acts freely on $\mu^{-1}(\tau)$, where
$$
     \tau:=\sum_{\nu=1}^n\w_\nu,
$$
and that the symplectic quotient 
$$
     \bar M :=  \C^n\dslash T(\tau) := \mu^{-1}(\tau)/T
$$
has minimal Chern number 
$
     N := \max\left\{m\in\Z\,|\,
     \tau/m\in\Lambda^*\right\}
     \ge 2.
$
These conditions guarantee that the symplectic
quotient $\bar M$ is a monotone toric manifold.
Combining Theorem~\ref{thm:invariants} with the 
results of~\cite{GS} one can compute the genus zero
Gromov--Witten invariants of the symplectic quotient. 
More precisely, denote by 
$
    S^*(\tt^*)\to H^*(\bar M;\R):\alpha\mapsto\bar\alpha
$
the Kirwan homomorphism. Consider the dual homomorphism 
$H_2(\bar M;\Z)\to\Lambda$ in degree two.  This homomorphism 
is injective.  We denote its image by $\Lambda(\tau)$
and the inverse map by 
$
     \Lambda(\tau)\to H_2(\bar M;\Z):\lambda\mapsto\bar\lambda.
$ 
Given $\lambda\in\Lambda(\tau)$ we denote by 
$\GW^{\bar M}_{\bar\lambda}$ the genus zero
Gromov--Witten invariant of $\bar M$ with fixed marked 
points in the homology class $\bar\lambda$.  
In~\cite[Theorem~A]{GS} it is proved that,
for every $\lambda\in\Lambda(\tau)$ and every 
$n$-tuple $\ell=(\ell_1,\dots,\ell_n)$ of nonnegative integers, 
\begin{equation}\label{eq:GS}
     \Phi^{\rho,\tau}_\lambda(\w^\ell)
     = \GW^{\bar M}_{\bar\lambda}
     (\bar\w_1,\dots,\bar\w_1,\dots,
      \bar\w_n,\dots,\bar\w_n),
\end{equation}
where each argument $\bar\w_\nu$ occurs $\ell_\nu$ times. 
Thus Theorem~\ref{thm:invariants} allows us to compute the 
genus zero Gromov--Witten invariants of tuples of cohomology 
classes of degree two. This can be used to compute the quantum 
cohomology ring of the symplectic quotient.  The statement 
of the theorem requires some preparation.  

The {\bf chamber} $C(\tau)$ 
is defined as the 
component of the set of regular values of $\mu$ that contains $\tau$.
The {\bf effective cone} $\Lambda_\eff(\tau)\subset\Lambda(\tau)$
is defined as the set of lattice vectors $\lambda\in\Lambda(\tau)$
that satisfy $\inner{\tau'}{\lambda}\ge 0$ for every 
$\tau'\in C(\tau)$.  

Let $\Rr$ be any graded commutative algebra (over the reals) with unit
which is equipped with a homomorphism
$$
     \Lambda_\eff(\tau)\to\Rr:\lambda\mapsto q^\lambda
$$
from the additive semigroup $\Lambda_\eff(\tau)$ to the 
multiplicative semigroup $\Rr$ such that 
$
     \deg(q^\lambda) = 2\inner{\tau}{\lambda}. 
$
Given such a graded algebra define the quantum cohomology 
ring $\QH^*(\bar M;\Rr)$ as the tensor product 
$$
     \QH^*(\bar M;\Rr) := H^*(\bar M;\R)\otimes\Rr
$$
(of vector spaces over the reals).  Thus an element of 
$\QH^*(\bar M;\Rr)$ is a finite sum   
$
     \bar\alpha = \sum_{r\in\Rr}
     \bar\alpha_r r
$
where $\bar\alpha_r\in H^*(\bar M;\R)$. 
The ring structure is defined by 
$$
     \bar\alpha'*\bar\alpha''
     := \sum_i\sum_{\lambda\in\Lambda_\eff(\tau)}
        \sum_{r',r''}
        \GW^{\bar M}_{\bar\lambda}
        (\bar\alpha'_{r'},\bar\alpha''_{r''},\bar e_i)
        \bar e_i^* r'r''q^\lambda,
$$
where the $\bar e_i$ form a basis of $H^*(\bar M;\R)$ and the 
$\bar e_i^*$ denote the dual basis with respect to the cup product
pairing.  

\begin{theorem}[Quantum Cohomology]\label{thm:batyrev}
Let $\bar M=\C^n\dslash T(\tau)$ 
be a (nonempty) monotone toric manifold with minimal
Chern number $N\geq 2$. Then the ring homomorphism
\begin{equation}\label{eq:uq}
      \Rr[u_1,\dots,u_n]\to\QH^*(\bar M;\Rr):
      ru^\ell\mapsto r\bar\w^{*\ell}
\end{equation}
induces an isomorphism 
$
      \QH^*(\bar M;\Rr)\cong \Rr[u_1,\dots,u_n]/\Jj,
$
where the ideal $\Jj\subset\Rr[u_1,\dots,u_n]$
is generated by the relations 
\begin{eqnarray*}
      \sum_{\nu=1}^n\eta_\nu\w_\nu=0
&\IMP &
      \sum_{\nu=1}^n\eta_\nu u_\nu=0,  \\
      \tau\notin C(\{1,\dots,n\}\setminus\{\nu\})
&\IMP &
      u_\nu=0, \\
     \lambda\in\Lambda_\eff(\tau),\;\;
     d_\nu^\pm:=\max\{\pm\inner{\w_\nu}{\lambda},0\}
&\IMP &
     \prod_\nu u_\nu^{d_\nu^+}
     = q^\lambda\prod_{\nu}u_\nu^{d_\nu^-}.
\end{eqnarray*}
\end{theorem}

The ring $\Rr[u_1,\dots,u_n]/\Jj$ was introduced by Batyrev~\cite{BATYREV}.
It also appeared in Givental's work on mirror 
symmetry~\cite{GIVENTAL} for the monotone case.  
Examples of Spielberg~\cite{SPIELBERG1,SPIELBERG2} show that, 
in the nonmonotone case, the kernel of the
homomorphism~(\ref{eq:uq}) is not necessarily equal to $\Jj$. 
For special cases the isomorphism 
$\Rr[u_1,\dots,u_n]/\Jj\to\QH^*(\bar M;\Rr)$
was established in~\cite{QR,SPIELBERG1}. 
The reason for our hypothesis $N\ge 2$ lies in the 
identity~(\ref{eq:GS}) which, in general, does not continue to hold
in the case $N=1$ (the degrees of all the classes must be less
than twice the minimal Chern number). 

In Section~\ref{sec:vortex} we explain some background 
from~\cite{CGMS} about the symplectic vortex equations.  
The wall crossing formula (for arbitrary genus) 
is restated in Section~\ref{sec:wall} and proved 
in Section~\ref{sec:proof}.  We prove Theorem~\ref{thm:invariants}  
in Section~\ref{sec:compute} and Theorem~\ref{thm:batyrev} in 
Section~\ref{sec:qcoh}.


\section{The symplectic vortex equations}\label{sec:vortex}

Fix a compact Riemann surface $(\Sigma,j_\Sigma,\dvol_\Sigma)$,
a principal $T$-bundle $P\to\Sigma$, and an invariant inner 
product on $\tt$. The characteristic vector
of $P$ will be denoted by
$$
\lambda(P) := \int_\Sigma F_A \in \Lambda.
$$
Here $A\in\Aa=\Aa(P)$ is a connection on $P$ and $\lambda(P)$
is independent of the choice of the connection.
The symplectic vortex equations (at a parameter $\tau\in\tt^*$)
have the form
\begin{equation}\label{eq:vortex}
\bar\p_Au_\nu=0,\qquad
*F_A+\pi\sum_{\nu=1}^n\left|u_\nu\right|^2\w_\nu
= \frac{*_\tt\lambda(P)}{\Vol(\Sigma)} + \tau,
\end{equation}
where $u_\nu:P\to\C$ is equivariant with respect to
$\rho_\nu$, $*:\Om^2(\Sigma,\tt)\to\Om^0(\Sigma,\tt^*)$
denotes the Hodge $*$-operator determined by the volume form 
on $\Sigma$ and the inner product on $\tt$,
and $*_\tt:\tt\to\tt^*$ denotes the isomorphism induced
by the inner product.  The gauge group
$
\Gg:=\Cinf(\Sigma,T)
$
acts on the space of solutions of~(\ref{eq:vortex})
with finite isotropy if and only if $\tau$ is a
regular value of $\mu$. Moreover, the moduli space
$$
\Mm(\tau):=\left\{(A,u)\,|\,u\mbox{ and }A
\mbox{ satisfy }(\ref{eq:vortex})\right\}/\Gg_0
$$
of based gauge equivalence classes of solutions of~(\ref{eq:vortex})
is compact (see~\cite{CGMS}). Here we fix a point $z_0\in\Sigma$ and
denote the based gauge group by
$$
\Gg_0 := \left\{g\in\Gg \,|\,g(z_0)=\one\right\}.
$$
Think of this moduli space as a subset
of the space
$$
\Bb := \frac{\Aa(P) \times\Cinf_T(P,\C^n)}{\Gg_0},
$$
where $g\in\Gg_0$ acts by $g^*(A,u):=(A+g^{-1}dg,\rho(g)^{-1}u)$.
The group $T$ (of constant gauge transformations) acts
contravariantly on $\Bb$ and we identify
$$
\Bb\times_T\ET \cong \frac{\Aa(P) \times\Cinf_T(P,\C^n)\times ET}
{\Gg}
$$
where $\Gg$ acts by
$g^*(A,u,e):=(A+g^{-1}dg,\rho(g)^{-1}u,g(z_0)^{-1}e)$.
The invariants introduced in~\cite{CGMS} are obtained by
integrating equivariant cohomology classes of $\Bb$
over $\Mm(\tau)/T$. Here integration is to be understood
in terms of evaluating the equivariant Euler class
of an associated $T$-moduli problem on the equivariant
cohomology class in question (see~\cite{CMS}).

It is useful to think of~$u_\nu$ as a section of the
complex line bundle
\begin{equation}\label{eq:Lnu}
L_\nu := P\times_{\rho_\nu}\C \to \Sigma,
\end{equation}
where the equivalence relation on $P\times\C$ is 
$[p,\zeta]\equiv[pg,\rho_\nu(g)^{-1}\zeta]$ for $g\in T$.
This bundle has degree
$$
d_\nu := \frac{i}{2\pi}\dot\rho_\nu(\lambda(P))
= \inner{\w_\nu}{\lambda(P)}
$$
and the moduli space $\Mm(\tau)/T$ has virtual dimension
$$
{\rm dim}_\R\Mm(\tau)/T
= (n-\dim\,T)(2-2g) + 2\sum_{\nu=1}^nd_\nu =: 2m,
$$
where $g$ is the genus of $\Sigma$.
Integration over the moduli space gives rise to
a homomorphism 
$$
     \Phi^{\rho,\tau}_{\lambda,g}:
     S^*(\tt^*)\otimes H^*(\Aa/\Gg_0)\to\R.
$$
Formally, this homomorphism can be written in the form 
\begin{equation}\label{eq:inv}
\Phi^{\rho,\tau}_{\lambda,g}(\alpha)
:= \int_{\Mm(\tau)/T}\pi^*\alpha
\end{equation}
for 
$
     \alpha\in S^*(\tt^*)\otimes H^*(\Aa/\Gg_0)\cong 
     H^*_T(\Aa/\Gg_0).
$ 
Here $\pi:\Bb\to\Aa/\Gg_0$ is the obvious projection
and $\pi^*:H^*_T(\Aa/\Gg_0)\to H^*_T(\Bb)$ is the 
induced homomorphism on equivariant cohomology. 
The precise definition of the invariant involves the 
equivariant Euler class of the $T$-moduli problem 
associated to equation~(\ref{eq:vortex}) 
(see~\cite{CMS,CGMS}).  The invariant 
$\Phi^{\rho,\tau}_{\lambda,g}(\alpha)$ 
can only be nonzero when $\alpha$ has degree $2m$ and 
it takes rational values on cohomology classes that 
take the integer lattice in $\tt$ to the integral 
cohomology of $\Aa/\Gg_0$. 


\section{Wall crossing}\label{sec:wall}

In this section we formulate the wall crossing formula
for arbitrary genus at an element $\tau_0\in\tt^*$ of a 
wall of codimension one.  This means that there exists 
an index set $I\subset\{1,\dots,n\}$ satisfying the following 
conditions.
\begin{description}
\item[(i)]
The subspace 
$
      W_I:=\SPAN\{\w_\nu\,|\,\nu\in I\}
      \subset\tt^*
$
has dimension $k-1$ and $\w_\nu\notin W_I$ for every 
$\nu\notin I$. 
\item[(ii)]
$
     \tau_0\in\mu(\C^I),
$
where
$
     \C^I := \{x\in\C^n\,|\,x_\nu=0\mbox{ for }\nu\notin I\}.
$
\item[(iii)]
If $J\subset\{1,\dots,n\}$ is another index set
satisfying~(i) then $\tau_0\notin\mu(\C^J)$.
\end{description}
Note that under these conditions $\tau_0$ is a positive linear 
combination of precisely $k-1$ linearly independent
vectors from the set $\{\w_\nu\,|\,\nu\in I\}$.
Choose a vector $\tau_1\in\tt^*$ 
that is transverse to $\mu(\C^I)$ and let
$e_1\in\Lambda$ be the unique primitive lattice vector 
that satisfies $\inner{\tau_1}{e_1}>0$ and is orthogonal 
to the wall at $\tau_0$:
$$
     \inner{\w_\nu}{e_1}=0\qquad\mbox{ for }\nu\in I.
$$
Denote by $T_1\subset T$ the subtorus generated by $e_1$ 
and by $\tt_1$ its Lie algebra. Let $T_0:=T/T_1$ 
be the quotient torus and $\tt_0:=\tt/\tt_1$ be its Lie algebra. 
Then the action $\rho$ induces an action $\rho_0$ of $T_0$
on $\C^I$.  

The wall crossing number will be expressed
as an integral over the moduli space $\Mm_0$ of based
gauge equivalence classes of solutions
$(A,\{u_\nu\}_{\nu\in I})$ of the equations
\begin{equation}\label{eq:v0}
\bar\p_Au_\nu=0 \quad (\nu\in I),\qquad
*F_A+\pi\sum_{\nu\in I}\left|u_\nu\right|^2\w_\nu
= \frac{*_\tt\lambda(P)}{\Vol(\Sigma)} + \tau_0.
\end{equation}
We shall view this as a $T_0$-moduli problem.
Indeed, the subgroup $T_1\subset T$ acts trivially 
on $\Mm_0$. However, since $\tau_0$ is a regular value 
of $\mu|_{\C^I}$, the quotient group $T_0=T/T_1$
acts on $\Mm_0$ with finite isotropy.

It is interesting to compare $\Mm_0$ with the
moduli space $\Mm(\C^I,P_0,\tau_0)$ of based gauge equivalence
classes of solutions of~(\ref{eq:vortex}) with
$\C^n$, $T$, and $P$ replaced $\C^I$, $T_0:=T/T_1$, 
and $P_0:=P/T_1$. There is a natural projection
\begin{equation}\label{eq:fibration}
\Mm_0\TO\Mm(\C^I,P_0,\tau_0):
[A,\{u_\nu\}_{\nu\in I}]\mapsto[A_0,\{u_\nu\}_{\nu\in I}],
\end{equation}
where $A_0:=\Pi_0A\in\Om^1(P,\tt_0)$ can be thought of
as a connection on $P_0$.  Here $\Pi_0:\tt\to\tt_0$ denotes the
canonical projection.  If $(A,\{u_\nu\}_{\nu\in I})$
satisfies~(\ref{eq:v0}) then the tuple 
$(A_0,\{u_\nu\}_{\nu\in I})$ satisfies the equations 
$$
\bar\p_{A_0}u_\nu=0 \quad (\nu\in I),\qquad
*_0F_{A_0}+\pi\sum_{\nu\in I}\left|u_\nu\right|^2\w_\nu
= \frac{*_{\tt_0}\lambda(P_0)}{\Vol(\Sigma)} + \tau_0.
$$
and hence belongs to the moduli space $\Mm(\C^I,P_0,\tau_0)$.
Here $*_{\tt_0}:\tt/\tt_1\to\tt_1^\perp$ is given by 
$
      [\xi]\mapsto \xi-|e_1|^{-2}\inner{\xi}{e_1}e_1
$
and $*_0:\Om^2(\Sigma,\tt/\tt_1)\to\Om^0(\Sigma,\tt_1^\perp)$
is induced by the Hodge $*$-operator on $\Sigma$ and $*_{\tt_0}$. 
The map~(\ref{eq:fibration}) defines a fibration whose fiber
can be described as follows.  

Choose a complement of $\tt_1$ in $\tt$ and denote the 
resulting projection by $\Pi_1:\tt\to\tt_1$.  
Define $\Pic_1$ as the space of based gauge equivalence 
classes of real valued $T$-invariant
$1$-forms $A_1\in\Om^1(P)$ that satisfy
$A_1(p e_1)=1$ for every $p\in P$ and
\begin{equation}\label{eq:vI}
*F_{A_1} = \frac{\Pi_1\lambda(P)}{\Vol(\Sigma)}.
\end{equation}
Here the based gauge group is
$
\Gg_1:=\left\{g_1:\Sigma\to T_1\,|\,g_1(z_0)=\one\right\}
$
and it acts by $({g_1}^*A_1)e_1:=A_1e_1+{g_1}^{-1}dg_1$.
Note that $\Pic_1$ is a $2g$-torus. It is the fibre
in~(\ref{eq:fibration}) because the subgroup
$\Gg_1\subset\Gg_0$ acts trivially on $A_0=\Pi_0A$ 
and $u_\nu$ for every $\nu\in I$. We emphasize that the
fibration~(\ref{eq:fibration}) need not be a product.
The reason is that the map $A\mapsto A_1:=\Pi_1A$
will not, in general, be gauge invariant.
We are now in a position to state the wall crossing 
formula for arbitrary genus.  Think of 
$\Mm_0$ as a subset of the space 
$\Bb_0:=\Aa(P)\times\Cinf_T(P,\C^I)/\Gg_0$
and denote by 
$
     \pi_0^*:S^*(\tt_0^*)\otimes H^*(\Aa/\Gg_0)
     \to H^*_{T_0}(\Bb_0)
$ 
the homomorphism on equivariant cohomology induced by 
the $T_0$-invariant projection $\pi_0:\Bb_0\to\Aa/\Gg_0$. 

\begin{theorem}[Wall Crossing]\label{thm:WCN}
Let $\alpha\in S^*(\tt^*)\otimes H^*(\Aa/\Gg_0)$, 
$\lambda\in\Lambda$, and $d_\nu:=\inner{\w_\nu}{\lambda}$.  
Then, for every sufficiently small positive number $\eps$, 
we have
\begin{equation}\label{eq:WCN}
\Phi^{\rho,\tau_0+\eps\tau_1}_{\lambda,g}(\alpha)
- \Phi^{\rho,\tau_0-\eps\tau_1}_{\lambda,g}(\alpha)
= \int_{\Mm_0/T_0}\pi_0^*\alpha_0,
\end{equation}
where $\alpha_0\in S^*(\tt_0^*)\otimes H^*(\Aa/\Gg_0)$ is 
the polynomial map defined by
$$
\alpha_0(\xi) := \frac{1}{2\pi i}\oint
\frac{\alpha(\xi+ze_1)}
{\prod_{\nu\notin I}\inner{\w_\nu}{\xi+ze_1}^{d_\nu+1-g}}
\exp\left(\sum_{\nu\notin I}\frac{\Om_\nu}
{\inner{\w_\nu}{\xi+ze_1}}\right)\,dz.
$$
\end{theorem}

Here $\Om_\nu$ is the closed $2$-form on $\Aa/\Gg_0$ defined by
$$
\Om_\nu := \sum_{j,j'=1}^k \Om_{jj'}\w_{\nu j}\w_{\nu j'}.
$$
Here we have chosen a basis $e_1,\dots,e_k$ of $\tt$ and a
symplectic basis $\alpha_1,\dots,\alpha_{2g}$ of $H^1(\Sigma;\Z)$.
These bases induce a basis $\tau_{ij}$ of $H^1(\Aa/\Gg_0;\Z)$
and $\Om_{jj'}$ and $\w_{\nu j}$ are defined by
$$
\Om_{jj'} := \sum_{i=1}^g\tau_{ij}\wedge\tau_{i+g,j'},\qquad
\w_{\nu j} := \inner{\w_\nu}{e_j}.
$$
Note that $\Om_{jj'}$ is independent of the choice of 
the $\alpha_i$ and $\Om_\nu$ is independent of the choice 
of both bases. For each $\xi$, the integral in the definition
of $\alpha_0(\xi)$ is over a circle in the complex plane
enclosing all the poles of the integrand.
The integral in~(\ref{eq:WCN}) is understood as the 
evaluation of the equi\-variant Euler class of the $T_0$-moduli
problem associated to~(\ref{eq:v0}) on the class $\pi_0^*\alpha_0$ 
(see~\cite{CMS}). 

\begin{remark}[Residues]\label{rmk:residue}\rm
Consider a rational function $f:\C\to\C$ with poles $p_1,\dots,p_n$. 
It induces a meromorphic 1-form $f\,dz$. 
Let $\oint f\,dz$ be the integral of $f\,dz$ over a closed
curve in $\C$ around all the poles of $f$. 
By the residue theorem,
$$
     \frac{1}{2\pi i}\oint f\,dz = -\Res_\infty(f\,dz)
     = \sum_{j=1}^n\Res_{p_j}(f\,dz). 
$$
Note that the 1-form $f\,dz$ and hence the residue at infinity do not 
change if the complex coordinate is shifted by $z\mapsto z+c$. 
If we expand $f$ as a Laurent series
$$
     f(z) = \sum_{k=-\infty}^{k_0} a_k z^k
$$
in $z^{-1}$ that converges near infinity then 
the residue at infinity is minus the coefficient of $z^{-1}$,
i.e. $\Res_\infty(f)=-a_{-1}$.
\end{remark}

\begin{example}\label{ex:wps}\rm
Consider the action of the 1-torus $T=\R/\Z$ on $\C^n$ with positive
integer weights 
$\w_\nu=\ell_\nu\in\Lambda^*\cong\Z$. 
Let $c\in\Lambda^*$
be the standard generator $c(\xi)=\xi$ and pick a homology class
$\lambda=d\in\Lambda\cong\Z$. 
Assume
$$
     m := \sum_{\nu=1}^n(d\ell_\nu+1-g)+g-1 \ge 0.
$$
We compute the invariant $\Phi^{\rho,\tau}_{d,g}(c^m)$
in the nonempty chamber by wall crossing from
the empty chamber. 
Here $I=\emptyset$, $T_1=T=\R/\Z$,
$e_1=1$, and $T_0=\{\one\}$.  Then 
$\Mm_0=\Pic$ is a $2g$-torus and 
$\Om_\nu={\ell_\nu}^2\Om$, where $\Om$ is
the standard symplectic form on $\Pic$. 
It satisfies
$$
\frac{1}{g!}\int_\Pic\Om^g = 1.
$$
The integrand in Theorem~\ref{thm:WCN} is given by
$$
\alpha_0 = \frac{1}{2\pi i}\oint
\frac{z^m}{\prod_\nu (\ell_\nu z)^{d\ell_\nu+1-g}}
\exp\left(\frac{(\ell_1+\dots+\ell_n)\Om}{z}\right)
= \frac{(\ell_1+\cdots+\ell_n)^g}
{\prod_{\nu=1}^n\ell_\nu^{d\ell_\nu+1-g}}
\frac{\Om^g}{g!}.
$$
Integrating this class over $\Pic$ yields the formula
from~\cite{CGMS}: 
$$
\Phi^{\rho,\tau}_{\lambda,g}(c^m)
= \frac{(\ell_1+\cdots+\ell_n)^g}
{\prod_{\nu=1}^n\ell_\nu^{d\ell_\nu+1-g}}.
$$
\end{example}


\section{Proof of the wall crossing formula}\label{sec:proof}


\subsection{A cobordism argument}

We introduce the gauge invariant differential equations
\begin{equation}\label{eq:vs}
\begin{gathered}
\bar\p_{A}u_\nu=0 \quad (\nu=1,\dots,n), \cr
*F_A+\pi\sum_{\nu\in I}\left|u_\nu\right|^2\w_\nu
= \frac{*_\tt\lambda(P)}{\Vol(\Sigma)} + \tau_0,\qquad
\sum_{\nu\notin I}\left\|u_\nu\right\|^2=1.
\end{gathered}
\end{equation}
Denote the moduli space of based gauge equivalence classes
of solutions of~(\ref{eq:vs})~by
$$
\Ss_0 := \left\{
(A,u)\,|\,A\mbox{ and }u\mbox{ satisfy }(\ref{eq:vs})
\right\}/\Gg_0.
$$
Note that there is a $T$-equivariant projection 
$\Ss_0\to\Mm_0$ whose fiber is the unit sphere in the 
kernel of the Cauchy--Riemann operator in the variables
$u_\nu$ for $\nu\notin I$.  

\begin{proposition}\label{prop:wcn}
The wall crossing number can be expressed in the form
$$
\Phi^{\rho,\tau_0+\eps\tau_1}_{\lambda,g}(\alpha)
- \Phi^{\rho,\tau_0-\eps\tau_1}_{\lambda,g}(\alpha)
= -\int_{\Ss_0/T}\pi^*\alpha,
$$
where the orientation of $\Ss_0/T$ is determined by $e_1$
as in Remark~\ref{rmk:ORIENT}.
\end{proposition}

\begin{remark}\label{rmk:ORIENT}\rm
The moduli space $\Ss_0/T$ is oriented as follows. Denote
$$
    \Xx:=\Om^1(\Sigma,\tt)\oplus\bigoplus_{\nu=1}^n\Om^0(\Sigma,L_\nu),\qquad
    \Yy:=\Om^0(\Sigma,\tt^*\otimes\C)
    \oplus\bigoplus_{\nu=1}^n\Om^{0,1}(\Sigma,L_\nu),
$$
where $L_\nu\to\Sigma$ is the line bundle~(\ref{eq:Lnu}).
Define the operator $\Dd_0:\Xx\to\Yy$ by 
$$
    \Dd_0(\alpha,\hat u)
    := \left(\begin{array}{c}
    *d\alpha+2\pi\sum_{\nu\in I}\inner{u_\nu}{\hat u_\nu}\w_\nu \\
    d^*\alpha+2\pi\sum_{\nu\in I}\inner{iu_\nu}{\hat u_\nu}\w_\nu \\
    \bar\p_A\hat u_\nu + \rho_\nu(\alpha)^{0,1}u_\nu
    \end{array}\right).
$$
The first component is the real part in 
$\Om^0(\Sigma,\tt^*\otimes\C)$ 
(it is the linearization 
of the second equation in~(\ref{eq:vs})) 
and the second component is the imaginary part (it corresponds 
to the local slice condition for action of $\Gg/T_1$). 
The third component is the linearization of the first equation
in~(\ref{eq:vs}).  The operator $\Dd_0$ is complex linear,
where the complex structure on $\Om^1(\Sigma,\tt)$ is given by 
the Hodge $*$-operator $\alpha\mapsto *\alpha=-\alpha\circ J_\Sigma$. 
Define the linear functional $\Psi:\Xx\to\C$ by 
$$
     \Psi(\alpha,\hat u)
     := \frac{2\pi}{\inner{\tau_1}{e_1}\Vol(\Sigma)}
        \sum_{\nu\notin I}\inner{\w_\nu}{e_1}
        \int_\Sigma\left(\inner{u_\nu}{\hat u_\nu}
        + i\inner{iu_\nu}{\hat u_\nu}\right)\,\dvol_\Sigma.
$$
The imaginary part of $\Psi$ corresponds to the local slice for
the $T_1$-action. 

Note that the cokernel of $\Dd_0$ always contains the constant
functions in 
$\Om^0(\Sigma,\tt_1\otimes\C)$.  
We assume for 
simplicity that the cokernel is equal to this space and that
$$
    \Psi(0,\{u_\nu\}_{\nu\notin I})
    = \frac{2\pi}{\inner{\tau_1}{e_1}\Vol(\Sigma)}
      \sum_{\nu\notin I}\left\|u_\nu\right\|^2_{L^2}\inner{\w_\nu}{e_1}
    \ne 0.
$$
Then $\Ss_0$ is a smooth manifold near $[A,u]$
and the tangent space of $\Ss_0/T$ is 
$$
    T_{[A,u]}\Ss_0/T
    = \left\{(\alpha,\hat u)\,|\,
      \Dd_0(\alpha,\hat u)=0,\,
      \sum_{\nu\notin I}\inner{u_\nu}{\hat u_\nu}=0,\,
      \IM\,\Psi(\alpha,\hat u)=0\right\}.
$$
A basis $v_1,\dots,v_{2m}$ of the tangent space is called
{\bf positively oriented} if the vectors
$v_1,\dots,v_{2m}$, $w_0:=(0,\{u_\nu\}_{\nu\notin I})$,
$w_1:=(0,\{-2\pi i\inner{\w_\nu}{e_1}u_\nu\}_{\nu\notin I})$
form a positive basis of the complex vector space $\ker\,\Dd_0$. 
Note that the orientation depends on $e_1$ and that 
$m=(n-\dim T)(1-g)+\sum_{\nu=1}^nd_\nu$. 
\end{remark}

\begin{proof}[Proof of Proposition~\ref{prop:wcn}.]
Denote $\tau_t:=\tau_0+t\tau_1$ and consider the moduli space
$$
\Ww := \left\{(t,A,u)\,|\,-\eps\le t\le\eps,\,
(\ref{eq:vortex})\mbox{ holds with }\tau=\tau_t\right\}/\Gg_0.
$$
This space has boundary
$$
\p\Ww = \bigl(\{-\eps\}\times\Mm(\tau_{-\eps})\bigr)
\cup \bigl(\{\eps\}\times\Mm(\tau_\eps)\bigr).
$$
Now the group $T$ of constant gauge transformations
does not act with finite isotropy on $\Ww$.
Call an element of $\Ww$ {\bf singular}
if its isotropy subgroup has positive dimension.
Then the singular subset of $\Ww$ is precisely the
moduli space $\Mm_0$ introduced in the previous section.
Cutting out a neighbourhood of this singular set
we obtain the subset $\Ww_\delta\subset\Ww$ of all
(based gauge equivalence classes of)
triples $[t,A,u]\in\Ww$ that satisfy
$
\sum_{\nu\notin I}\left\|u_\nu\right\|^2\ge\delta.
$
The boundary of this set is
$$
\p\Ww_\delta
= \bigl(\{-\eps\}\times\Mm(\tau_{-\eps})\bigr)\cup
\bigl(\{\eps\}\times\Mm(\eps)\bigr)\cup\Mm_\delta,
$$
where the third boundary component $\Mm_\delta$ is the moduli
space of based gauge equivalence classes
of solutions of the equations
\begin{equation}\label{eq:vdelta}
\begin{gathered}
\bar\p_{A}u_\nu=0 \quad (\nu=1,\dots,n),\cr
*F_A+\pi\sum_{\nu=1}^n\left|u_\nu\right|^2\w_\nu
= \frac{*_\tt\lambda(P)}{\Vol(\Sigma)} + \tau_t,\qquad
\sum_{\nu\notin I}\left\|u_\nu\right\|^2=\delta.
\end{gathered}
\end{equation}
Here $\left\|\cdot\right\|$ denotes the $L^2$-norm.
If $\delta>0$ is sufficiently small then $\Mm_\delta$ does not
intersect the boundary of $\Ww$. We point out that
\begin{equation}\label{eq:t}
t = \frac{\pi}{\inner{\tau_1}{e_1}\Vol(\Sigma)}
\sum_{\nu\notin I}\left\|u_\nu\right\|_{L^2}^2\inner{\w_\nu}{e_1}.
\end{equation}
To see this take the inner product of the second equation
in~(\ref{eq:vdelta}) with $e_1$ and integrate over $\Sigma$.
This shows that the parameter $t$ is determined by $u$
and can therefore be removed in the definition of $\Mm_\delta$.
If we orient $\Mm_\delta$ as the boundary of $\Ww_\delta$ then
the wall crossing number is given by 
$$
\Phi^{\rho,\tau_0+\eps\tau_1}_{\lambda,g}(\alpha)
- \Phi^{\rho,\tau_0-\eps\tau_1}_{\lambda,g}(\alpha)
= - \int_{\Mm_\delta/T}\pi^*\alpha.
$$
It remains to prove that the integrals over 
$\Mm_\delta/T$ 
and $\Ss_0/T$ agree. 

We prove this in two steps. The first step shows that the 
($T$-moduli problems associated to) $\Mm_\delta$ and $\Ss_0$
are cobordant and the second step is to compare the orientations. 
Fix a real parameter $s\in[0,1]$ and replace the second
equation in~(\ref{eq:vdelta}) by
\begin{equation}\label{eq:homotopy}
*F_A + \pi\sum_{\nu\in I}\left|u_\nu\right|^2\w_\nu
- \frac{*_\tt\lambda(P)}{\Vol(\Sigma)} - \tau_0
= s\left(
  t\tau_1-\pi\sum_{\nu\notin I}\left|u_\nu\right|^2\w_\nu
  \right).
\end{equation}
For $s=1$ this is equivalent to the second equation
in~(\ref{eq:vdelta}) and for $s=0$ to the second equation
in~(\ref{eq:vs}). Thus the moduli spaces $\Mm_\delta$ and
$\Ss_0$ are equivariantly cobordant and the isotropy subgroups
are finite for every $s$. 

Next we compare the orientation of $\Ss_0/T$ 
with the boundary orientation of $\Mm_\delta/T$. 
Let $[t,A,u]\in\Ww_\delta$ and define 
the operator $\Dd:\Xx\to\Yy$ by 
$$
    \Dd(\alpha,\hat u)
    := \left(\begin{array}{c}
    *d\alpha+2\pi\sum_{\nu=1}^n\inner{u_\nu}{\hat u_\nu}\w_\nu \\
    d^*\alpha+2\pi\sum_{\nu=1}^n\inner{iu_\nu}{\hat u_\nu}\w_\nu \\
    \bar\p_A\hat u_\nu + \rho_\nu(\alpha)^{0,1}u_\nu
    \end{array}\right).
$$
This operator is complex linear. Assume for simplicity 
that $\Dd$ is surjective.  Then $\Ww_\delta$ is a 
manifold near $[t,A,u]$ and the tangent space of 
$\Ww_\delta/T$ is 
$$
    T_{[t,A,u]}\Ww_\delta/T
    = \left\{(\hat t,\alpha,\hat u)\,|\,
      \Dd(\alpha,\hat u)=(\hat t\tau_1,0,0)\right\}
$$
The condition $\Dd(\alpha,\hat u)=(\hat t\tau_1,0,0)$ implies
$
    \Psi(\alpha,\hat u) = \hat t.
$
(Take the inner product of the first two components of 
$\Dd(\alpha,\hat u)$ with $e_1$ and integrate over $\Sigma$.)
Thus we may identify
$$
    T_{[t,A,u]}\Ww_\delta/T
    \cong \left\{(\alpha,\hat u)\,|\,
      \Dd(\alpha,\hat u)\in(\R\tau_1,0,0)\right\}.
$$
Let $v_1,\dots,v_{2m},w_0$ be a basis of the tangent space
with $\RE\Psi(w_0)>0$ and $v_\nu\in\ker\,\Dd$.
The basis is called {\bf positively oriented} if the vectors 
$v_1,\dots,v_{2m}$ form a positively oriented basis of 
the kernel of $\Dd$. 

Using the determinant bundle over the space of Fredholm 
operators we can carry the orientation of $\Ww_\delta$ 
over to a point $[A,u]\in\Ss_0$. 
The notion of an inward pointing vector also carries over.
Let $v_1,\dots,v_{2m}$ be a positive basis of $T_{[A,u]}\Ss_0/T$,
consider the inward pointing vector
$w_0:=(0,\{u_\nu\}_{\nu\notin I})$, and denote 
$w_1:=(0,\{-2\pi i\inner{\w_\nu}{e_1}u_\nu\}_{\nu\notin I})$.
Then, by Remark~\ref{rmk:ORIENT}, the vectors 
$v_1,\dots,v_{2m},w_0,w_1$ form a positive basis of 
$\ker\,\Dd_0$. Deform this basis continuously into a basis
$\hat v_1,\dots,\hat v_{2m},w_0,w_1$ of $\ker\,\Dd_0$
such that $\hat v_1,\dots,\hat v_{2m}$
form a basis of 
$\ker\,\Dd_0\cap\ker\,\Psi$. 
Note that $\IM\,\Psi(w_1)<0$. 
If $\RE\,\Psi(w_0)<0$ it follows that the vectors 
$\hat v_1,\dots,\hat v_{2m}$ form a positive basis
of 
$\ker\,\Dd_0\cap\ker\,\Psi$ 
and otherwise they form a negative basis. 
Thus, in both cases, the vectors 
$\hat v_1,\dots,\hat v_{2m},w_0$ represent
the negative orientation of the tangent space of $\Ww$. 
Since $w_0$ is inward pointing, it follows that the 
orientation of $\Ss_0/T$ agrees with the boundary orientation
of $\Mm_\delta/T$. Hence the integral of $\pi^*\alpha$
over $\Mm_\delta/T$ agrees with the integral over $\Ss_0/T$.

The assertion about the integrals can be understood 
verbatim if the moduli spaces are smooth.  
Otherwise they are understood in terms of the equivariant Euler 
classes of the associated $T$-moduli problems. 
The result then follows from the cobordism axiom 
for the Euler class in~\cite{CMS}.
\end{proof}

The moduli space $\Ss_0$ is the sphere
bundle in the kernel bundle of a family of
Cauchy--Riemann operators over $\Mm_0$.
In Section~\ref{sec:localize} we explain a general
equivariant localization formula for such kernel bundles.
The relevant index computation uses the
Atiyah--Singer index theorem for families
and will be carried out in Section~\ref{sec:index}.
The next section explains the necessary background 
about the equivariant Euler class.


\subsection{The equivariant Euler class}\label{sec:euler}

We begin with some recollections about
the equivariant Euler class (see~\cite{CMS} for details).
Let $X$ be a compact oriented smooth manifold,
$E\to X$ be an oriented real vector bundle of rank $k$,
and $\G$ be a compact Lie group which acts on $X$
and $E$ by orientation preserving diffeomorphisms such that
the projection is equivariant and the action is linear
on the fibres. We shall think of the action of $\G$
on $X$ and $E$ as a right action and denote it by
$(x,e)\mapsto(g^*x,g^*e)$ for $e\in E_x$.
The corresponding covariant action will be denoted
by $g_*x:=(g^{-1})^*x$ and the infinitesimal
(contravariant) action of $\xi\in\g:=\Lie(\G)$
by $\xi^*x\in T_xX$. An {\bf equivariant Thom form}
is a $d_\G$-closed equivariant differential form
$\tau_\G(E)\in\Om_\G^k(E)$
with compact support and fibre integral one.
The {\bf equivariant Euler class} $e_\G(E)\in H^k_\G(X)$ is the
cohomology class of the pullback of an equivariant Thom form under the
zero section. We will sometimes use the same notation for the Euler
class and a form representing it. The Thom class and the Euler class
are multiplicative under direct sum.

Now suppose that $E$ is a rank $n$ complex vector bundle and the
action of $\G$ is complex linear on the fibres. Then an explicit
representative of the equivariant Euler class can be constructed
as follows.  Fix a $\G$-invariant Hermitian metric on $E$
and let $P\to X$ denote the unitary frame bundle of $E$.
This bundle carries a right action of $\U(n)$
and a left action of $\G$, and these two action
commute. Let $X_\xi\in\Vect(P)$ denote the
infinitesimal (covariant) action of $\xi\in\g$.
More precisely a point $p\in P_x$ of the fibre over
$x\in X$ is a unitary vector space isomorphism
$p:\C^n\to E_x$ and the left action of $g\in\G$ is
given by $g_*p:\C^n\to E_{g_*x}$.
The vector field $X_\xi\in\Vect(P)$ is defined by
$$
     X_\xi(p) := \left.\frac{d}{dt}\right|_{t=0}
     \exp(t\xi)_*p
     \in T_pP.
$$

\begin{lemma}\label{le:euler-conn}
Let $A\in\Aa(P)\subset\Om^1(P,\u(n))$ be a
$\G$-invariant $\U(n)$-con\-nec\-tion form on $P$.
Then the $\G$-equivariant Euler class of a complex vector
bundle $E$ is represented by the $d_\G$-closed form
\begin{equation}\label{eq:euler-conn}
e_\G(E,\xi)
= \det\left(\frac{i}{2\pi}F_A +
  \frac{i}{2\pi}A(X_\xi)\right),
\end{equation}
where $F_A\in\Om^2(P,\u(n))$ denotes the curvature of~$A$.
\end{lemma}

\begin{proof}
The right hand side in~(\ref{eq:euler-conn})
is invariant and horizontal for the $\U(n)$-action
and thus descends to a $\G$-equivariant form on $X$. It is easy to
check that this form is $d_\G$-closed and hence represents an
equivariant cohomology class.

To prove~(\ref{eq:euler-conn}) we assume first that
$E=X\times\C^n$ is a trivial bundle and $\rho:\G\to\U(n)$
is a unitary representation of $\G$. The homomorphism
$\rho$ defines the covariant action of $\G$ on $E$
and so
$$
    g^*(x,z) = (g^*x,\rho(g)^{-1}z)
$$
for $x\in X$, $z\in\C^n$, and $g\in\G$. The frame bundle
of $E$ is the product bundle $P:=X\times\U(n)$ and
the formula
$$
    A_{x,u}(v,u\eta) := \eta
$$
for $v\in T_xX$, $u\in\U(n)$ and $\eta\in\u(n)$
defines a $\U(n)$-connection form $A\in\Om^1(P,\u(n))$.
This connection is $\G$-invariant and flat. For $\xi\in\g$
the vector field $X_\xi\in\Vect(P)$ is given by
$$
    X_\xi(x,u) = (0,\dot\rho(\xi)u)
$$
and so
$$
    A(X_\xi(x,u)) = u^{-1}\dot\rho(\xi)u,\qquad
    \det\left(\frac{i}{2\pi}A(X_\xi)\right)
    = \det\left(\frac{i}{2\pi}\dot\rho(\xi)\right).
$$
Now a Thom form on $E$ can be constructed as follows.
For $k=0,\dots,n$ let $\sigma_k:\u(n)\to\Om^{2n-2k}(\C^n)$
be a polynomial map of degree $k$.  It is shown
in~\cite[Lemma~5.5]{CMS} that these polynomials can be
chosen such that $\sigma_0\in\Om^{2n}(\C^n)$ is the standard
volume form,
$$
    \sigma_n(\eta) = \det(i\eta),
$$
and
$$
    \i(v_\eta)\sigma_k(\eta) = \lambda\wedge\sigma_{k+1}(\eta)
$$
for each $k$, where $\lambda\in\Om^1(\C^n)$
is the differential of the function $z\mapsto|z|^2/2$
and the vector field $v_\eta\in\Vect(\C^n)$ is defined by
$v_\eta(z):=\eta z$ for $\eta\in\u(n)$.
Now choose functions
$f_k:[0,\infty)\to[0,\infty)$ with compact support
such that $f_0(s)=0$ for $s\le\delta$ and $s\ge 1$
and
$$
    f_k'(s) + f_{k-1}(s) = 0,\qquad f_k(1)=0,
$$
and
$$
    \int_0^\infty s^kf_0(s)\,ds = 0,\qquad
    \int_0^\infty s^{n-1}f_0(s)\,ds
    = \frac{1}{2^{n-1}\Vol(S^{2n-1})}
$$
for $0\le k\le n-2$. Then
$$
     f_k(s)=\frac{1}{(k-1)!}\int_s^1(t-s)^{k-1}f_0(t)\,dt
$$
and hence $f_k(s)=0$ for $s<\delta$ and $k<n$ and
$$
    f_n(0)
    = \frac{1}{2^{n-1}(n-1)!\Vol(S^{2n-1})}
    = \frac{1}{(2\pi)^n}.
$$
Now a Thom form on $E=X\times\C^n$ is given by
$$
    \tau(\xi) = \sum_{k=1}^nf_k(|z|^2/2)\sigma_k(\dot\rho(\xi)).
$$
Its pullback under the zero section is given by
$$
    e_\G(E)
    = f_n(0)\sigma_n(\dot\rho(\xi))
    = \det\left(\frac{i}{2\pi}\dot\rho(\xi)\right).
$$
This proves the lemma in the case $E=X\times\C^n$.
For general $\G$-equivariant bundles $E\to X$ the result
follows from the {\it (Naturality)} axiom for the Euler
class and the fact that the pullback of $E$ under the projection
$P\to E$ is isomorphic to the trivial bundle $P\times\C^n$.
\end{proof}

\begin{remark}\rm
The formula of Lemma~\ref{le:euler-conn} can also
be expressed as follows. Let $\nabla$ be a $\G$-equivariant
Hermitian connection on $E$ and, for $\xi\in\g$,
denote by $\xi^\nabla\in\Om^0(X,\End(E))$
the covariant infinitesimal action
defined by
$$
    \xi^\nabla e:= \Nabla{t}\exp(t\xi)_*e|_{t=0}
$$
Then the Euler class is given by
$$
e_\G(E,\xi)
= \det\left(\frac{i}{2\pi}F^\nabla +
  \frac{i}{2\pi}\xi^\nabla\right).
$$
\end{remark}

\begin{example}\label{ex:euler-diag}\rm
Let $E\to X$ be a rank $n$ complex vector bundle.
Suppose $\G$ acts trivially on $X$ and
that the covariant action on the fibres is given
by a homomorphism $\rho:\G\to S^1$, given by
$$
     \rho(\exp(\xi))=e^{-2\pi i\inner{\w}{\xi}}, 
$$
where $\w\in\g^*$. Then, for every $\G$-invariant Hermitian
connection $\nabla$ on $E$, the endomorphism
$\xi^\nabla\in\Om^0(X,\End(E))$ is given by multiplication
with the imaginary number $\dot\rho(\xi)$.  Hence
$$
     e_\G(E,\xi) 
     = \det\left(
       \frac{i}{2\pi}F^\nabla +\frac{i}{2\pi}\dot\rho(\xi)
       \right)
     = \sum_{j=0}^n\inner{\w}{\xi}^{n-j}c_j(E).
$$
\end{example}

We wish to invert the equivariant Euler class. This requires an 
extension of the equivariant cohomology ring of $X$.

\smallskip
\noindent{\bf Standing assumption.}
{\it In the following $X$ is a smooth manifold, 
$\G$ is a compact Lie group acting on $X$,
and $T_1\subset\G$ is an oriented circle which is 
contained in the center of $\G$ and acts trivially on $X$.}
\smallskip

Denote the quotient group by $\G_0:=\G/T_1$.
Denote by $\e_1$ the positive integral generator 
of the Lie algebra $\tt_1:=\Lie(T_1)$.
Let $n$ be an integer. A {\bf $T_1$-rational $\G$-equivariant 
differential form of degree $n$ on $X$} 
is a Laurent series in $z^{-1}$ of the form 
$$
     \alpha(\xi,z) = \sum_{j\le n/2}\alpha_j(\xi)z^j,
$$
with coefficients $\alpha_j\in\Om^{n-2j}_\G(X)$,
that satisfies the following conditions.

\smallskip
\noindent{\bf (i)}
For every $\xi\in\g$ and every $x\in X$ the Laurent series 
$\sum_{j\le n/2}\alpha_j(\xi)_xz^j$
is a rational function on $\C$ with values 
in the complex vector space $\Lambda^*T_x^*X\otimes\C$. 

\smallskip
\noindent{\bf(ii)}
For every $t\in\R$ 
we have $\alpha(\xi,t+z)=\alpha(\xi+te_1,z)$.
Equivalently,
\begin{equation}\label{eq:alpha-1}
      \alpha_k(\xi+te_1) = \sum_{k\le j\le n/2}{j\choose k}
      \alpha_j(\xi)t^{j-k},\qquad k\ge 0,
\end{equation}
\begin{equation}\label{eq:alpha-2}
      \alpha_k(\xi+te_1) = \sum_{k\le j<0}{{-k-1}\choose{-j-1}}
      \alpha_j(\xi)(-t)^{j-k},\qquad k<0.
\end{equation}

Denote by $\Om^n_{\G,T_1}(X)$ the space of $T_1$-rational
$\G$-equivariant differential forms on $X$. This is a 
chain complex with respect to the usual equivariant
differential $d_\G\alpha(\xi):=d\alpha(\xi)+\i(X_\xi)\alpha(\xi)$.
The cohomology of this chain complex will be denoted by 
$H^*_{\G,T_1}(X)$.  

Let $\alpha=\sum_j\alpha_jz^j\in\Om^n_{\G,T_1}(X)$.
Then $\alpha_{-1}(\xi+te_1)=\alpha_{-1}(\xi)$. 
In other words, the coefficient of $z^{-1}$ 
descends to a $\G_0$-equivariant cohomology
class on $X$.  Minus this coefficient is called 
{\bf the residue at infinity of $\alpha$}
and will be denoted by 
$$
      \Res_\infty(\alpha)
      := -\alpha_{-1}(\xi)
      = -\frac{1}{2\pi i}\oint \alpha(\xi,z)\,dz
      \in \Om^{n+2}_{\G_0}(X).
$$
The residue at infinity descends to a homomorphism 
$$
     \Res_\infty:H^n_{\G,T_1}(X)\to H^{n+2}_{\G_0}(X).
$$

\begin{remark}\label{rmk:ratH3}\rm
There is an obvious inclusion $\Om^*_\G(X)\to\Om^*_{\G,T_1}(X)$
whose image is the subspace of polynomials
$\alpha=\sum_{0\le j\le n/2}\alpha_jz^j\in\Om^n_{\G,T_1}(X)$.  
Condition~(\ref{eq:alpha-1}) shows that any such form is 
uniquely determined by $\alpha_0\in\Om^n_\G(X)$ and vice versa. 
The inclusion $\Om^*_\G(X)\INTO\Om^*_{\G,T_1}(X)$
induces an inclusion in 
cohomology
$$
     H^*_\G(X)\INTO H^*_{\G,T_1}(X)
$$
whose left inverse is induced by the projection
$\alpha=\sum_j\alpha_jz^j\mapsto\alpha_0$. 
\end{remark}

Let $E\to X$ be a $\G$-equivariant complex vector bundle 
of rank $n_E$. The subgroup $T_1$ acts on $E$ with {\bf weight}
$$
     \w_E := \det\left(\frac{i}{2\pi}\dot\rho_x(e_1)\right)\in\Z.
$$
Here the homomorphism $\rho_x:T_1\to\Aut(E_x)$ 
denotes the action on the fiber over~$x$ and 
$\dot\rho_x:\tt_1\to\End(E_x)$ denotes the corresponding 
Lie algebra homomorphism.  The weight $\w_E$ is independent 
of $x$. Think of the equivariant Euler class as a polynomial map
$\g\to\Om^*(X)$. By Lemma~\ref{le:euler-conn}, 
the $\G$-equivariant Chern classes
$c_j(E)\in H^{2j}_\G(X)$ are the coefficients 
of $z^{n_E-j}$ in the polynomial 
$$
     e_\G(E,\xi+ze_1) =: \sum_{j=0}^{n_E}c_j(E,\xi)z^{n_E-j}.
$$
In particular,
$
     c_0(E,\xi)=\w_E.
$
If $\w_E\ne 0$ then the equivariant Euler class 
$e_\G\in H^{2n_E}_T(X)$ has a well
defined inverse $1/e_\G$ in the $T_1$-rational
$\G$-equivariant cohomology group $H^{-2n_E}_{\G,T_1}(X)$.
To see this, expand the rational function $z\mapsto 1/e_\G(\xi+ze_1)$
into a Laurent series in $z^{-1}$ which converges near infinity:
\begin{eqnarray*}
  \frac{1}{e_\G(E,\xi+ze_1)}
&= &
  \frac{1}{\w_Ez^{n_E}}
  \sum_{k=0}^\infty
  \left(-\sum_{j=1}^{n_E}\frac{c_j(E,\xi)}{\w_E}z^{-j}\right)^k \\
&=: &
  \sum_{i=0}^\infty s_i(E,\xi)z^{-n_E-i}.
\end{eqnarray*}
The coefficients $s_i(E)\in H^{2i}_\G(X)$ of this Laurent series 
are called the {\bf equi\-variant Segre classes} of $E$.
They are uniquely determined by the equation
\begin{equation}\label{eq:cs}
   \sum_{i+j=k}s_i(E,\xi)c_j(E,\xi)
   = \left\{\begin{array}{rl}
     1,&\mbox{if }k=0,\\
     0,&\mbox{if }k>0.
     \end{array}\right.
\end{equation}
In particular, the degree zero Segre class is 
$
     s_0(E,\xi) = 1/\w_E.
$
If $F\to X$ is another $\G$-equivariant complex vector bundle 
of rank $n_F$ with weight $\w_F$, then the quotient 
$$
     e_\G(F\ominus E) := \frac{e_\G(F)}{e_\G(E)}
     \in H^{2n_F-2n_E}_{\G,T_1}(X)
$$
depends only on the equivariant $K$-theory class 
$F\ominus E\in K_\G(X)$. It is only defined 
for equivariant $K$-theory classes $F\ominus E$ whose 
denominator $E$ has nonzero weight.


\subsection{Localization}\label{sec:localize}

Let $X$ be an orientable smooth manifold,
$\G=T$ be a torus acting on $X$, and 
$T_1\subset T$ be an oriented circle 
that acts trivially on $X$. 
We assume that the quotient group 
$T_0:=T/T_1$ acts on $X$ with finite isotropy. 
Denote the Lie algebras by $\tt:=\Lie(T)$, 
$\tt_1:=\Lie(T_1)$, and $\tt_0:=\tt/\tt_1:=\Lie(T_0)$,
let $\Lambda\subset\tt$ be the integer lattice,
and denote by $e_1\in\tt_1\cap\Lambda$
the positive generator of the sublattice.
Throughout we denote $m:=\dim\,X-\dim\,T_0$.

Let $\Ee\to X$ and $\Ff\to X$ be complex Hilbert space 
bundles on which $T$ acts complex linearly such that the 
projections are equivariant. Assume that $T_1$ acts with finite
isotropy outside the zero sections of $\Ee$ and $\Ff$.
Let
$$
\Dd_x:\Ee_x\to\Ff_x
$$
be a smooth family of $\G$-equivariant complex linear
Fredholm operators of complex numerical index
$$
\INDEX(\Dd) := \dim^\C\ker\,\Dd_x-\dim^\C\coker\,\Dd_x.
$$
Denote by
$$
\IND(\Dd) := \bigcup_{x\in X}\{x\}\times
\ker\,D_x\ominus\coker\,\Dd_x \in K_\G(X)
$$
the topological index of $\Dd$,
understood as a $\G$-equivariant $K$-theory class.
Consider the following $\G$-moduli problem.
The Hilbert manifold $\B$ is given by
$$
\B := \left\{(x,e)\,\Big|\,x\in X,\,
e\in\Ee_x,\,\left\|e\right\|^2=1
\right\},
$$
the Hilbert space bundle $\HH\to\B$ has fibre
$$
\HH_{x,e} := \Ff_x
$$
over $(x,e)\in\B$, and the section $\SS:\B\to\HH$ is given by
$$
\SS(x,e) := \Dd_xe.
$$
The zero set of this section is the {\it kernel manifold}
$$
\M := \left\{(x,e)\in\B\,|\,D_xe = 0\right\}.
$$
Denote by $\pi:\B\to X$ the obvious projection.
The equivariant $K$-theory class $\IND(\Dd)\in K_\G(X)$
has a nonzero weight (for the $T_1$-action) and hence
carries an equivariant Euler class
$$
e_\G(\IND(\Dd))\in H^*_{T,T_1}(X).
$$
in the $T_1$-rational $T$-equivariant cohomology of $X$.
The following theorem generalizes the localization formula
for circle actions in~\cite{CMS}. The assertion requires a
choice of orientations.

\begin{remark}\label{rmk:orient}\rm
Orientations of $X/T_0$ and $T_1$ determine
an orientation of the $T$-moduli problem
$(\B,\HH,\SS)$ as follows.  By choosing local
trivializations we may assume that $\Ee$ and $\Ff$
are (complex) Hilbert spaces equipped with a $T$-action
and so $\Dd$ is a $T$-equivariant smooth map
$X\to\Ll(\Ee,\Ff):x\mapsto\Dd_x$,
which assigns a (complex linear) Fredholm operator $\Dd_x$ 
to every $x\in X$. In this case the vertical differential
of $\SS$ at a point $(x,e)\in\M$ is an operator
$$
    D\SS(x,e):\left\{(\hat x,\hat e)\in T_xX\times\Ee\,|\,
    \inner{\hat e}{e}=0\right\}\to\Ff.
$$
It is given by
$$
    D\SS(x,e)(\hat x,\hat e) = \Dd_x\hat e + \dot\Dd(\hat x)e,
$$
where $\dot D(\hat x)e$ is defined as the derivative of
the path $\R\to\Ff:t\mapsto\Dd_{\exp_x(t\hat x)}e$ at $t=0$.
Now suppose that $\Dd_x$ is surjective.
Then a {\bf positive basis} of the kernel
of $D\SS(x,e)$ is defined as follows.
Pick a positive basis $\hat x_1,\dots,\hat x_m$
of $T_xX/\tt_0x$ and choose $\hat e_1,\dots,\hat e_m\in\Ee$
such that $\Dd_x\hat e_i + \dot\Dd(\hat x_i)e=0$
and $\inner{\hat e_i}{e}=0$ for $i=1,\dots,m$.  Next 
choose a positive basis $\hat e_{m+1},\dots,\hat e_{m+2n}$ of
the complex vector space $\ker\,\Dd_x$ such that
$\hat e_{m+2n-1}=e$ and $\hat e_{m+2n}$
is a positive tangent vector of the $T_1$-orbit
of $e$. Then the vectors $(\hat x_i,\hat e_i)$ for
$i=1,\dots,m$ and the vectors $(0,\hat e_j)$
for $j=m+1,\dots,m+2n-2$ are declared to be
a positive basis of $\ker\,D\SS(x,e)/\tt\cdot(x,e)$.
This definition of the orientation is independent
of the choices. If $\Dd_x$ is not surjective,
one can apply the same construction to the
kernel of a suitably augmented operator.
We emphasize that the orientation described here 
agrees with the convention of Remark~\ref{rmk:ORIENT}. 
\end{remark}

\begin{theorem}\label{thm:local}
Let $m:=\dim\,X-\dim\,T_0$ and $n:=\INDEX(\Dd)$.
Fix any orientation of $X/T_0$, 
let $T_1$ be oriented by $e_1$,
and orient $\M/T$ as in Remark~\ref{rmk:orient}. 
Then
\begin{equation}\label{eq:local}
    \int_{\M/T}\pi^*\alpha 
    = \int_{X/T_0}\Res_\infty
      \left(\frac{\alpha}{e_T(\IND(\Dd)}\right)
\end{equation}
for every $\alpha\in H_T^{m+2n-2}(X)$.
\end{theorem}

The integral on the left is understood as the
Euler class of the $T$-moduli problem $(\B,\HH,\SS)$
evaluated on $\pi^*\alpha$ (see~\cite{CMS}). 
The integrand on the right is the residue at infinity
of the $T_1$-rational $T$-equivariant cohomology 
class $\alpha/e_T(\IND(\Dd))\in H^{m-2}_{T,T_1}(X)$.
It is a $T_0$ equivariant cohomology class in $H^m_{T_0}(X)$
and can be integrated over $X/T_0$ because $T_0$ acts on $X$
with finite isotropy. 

\begin{remark}\label{rmk:local}\rm
Theorem~\ref{thm:local} continues to hold if we replace 
$X$ by a $T_0$-moduli problem $(\Bb_0,\Ee_0,\Ss_0)$
as in~\cite{CMS} and $\Ee$ and $\Ff$ by Hilbert space bundles 
over $\Bb_0$.  Then $\B$ is the unit 
sphere bundle in $\Ee$, $\HH_{b,e}=\Ee_{0b}\oplus\Ff_b$,
$\SS(b,e)=(\Ss_0(b),\Dd_be)$, and the right hand side 
of~(\ref{eq:local}) is understood in terms of the Euler
class of $(\Bb_0,\Ee_0,\Ss_0)$.  To prove this, choose
a finite dimensional reduction of $(\Bb_0,\Ee_0,\Ss_0)$
and note that~(\ref{eq:local}) continues to hold for 
noncompact manifolds $X$ and compactly supported 
$T$-equivariant differential forms $\alpha$.
\end{remark}

\begin{proof}[Proof of Theorem~\ref{thm:local}.]
The proof has three steps.

\medskip
\noindent{\bf Step~1.}
{\it We may assume without loss of generality
that $\Ee$ is finite dimensional and admits an equivariant
trivialization and that $\Ff=0$.
}

\medskip
\noindent
The reduction to the finite dimensional case
is proved as in~\cite[Theorem 11.1]{CMS}.
Hence assume $E=\Ee$ and $F=\Ff$ are finite dimensional.
By Proposition~\ref{prop:triv}, there exists a $T$-equivariant
complex vector bundle $E'\to X$ such that $E\oplus E'$ is
equivariantly isomorphic to $X\times V$ for some complex
$T$-representation~$V$.   Since $T$ is a torus
there exists a homomorphism $T\to S^1$ whose
restriction to $T_1$ has nonzero degree.
Multiplying the action of $T$ on $E'$ with a suitable
power of this homomorphism we may assume that
the action of $T_1$ on $E'$ has nonzero weight.
Now let $\B'\subset E\oplus E'$ be the unit sphere
bundle, $\HH'\to\B'$ be the pullback of $F\oplus E'$
under the projection $\pi':\B'\to X$,
and $\SS':\B'\to\HH'$ be given by
$$
    \SS'(x,e,e'):=(x,D_xe,e').
$$
Then the inclusion 
$
    \B\to\B':(x,e)\mapsto(x,e,0)
$
defines a morphism of $T$-moduli problems.
Hence, assuming the assertion for $E$ replaced by
$X\times V$ and $F$ replaced by the zero bundle,
we obtain
\begin{eqnarray*}
\int_{\M/T}\pi^*\alpha
&= &
\chi^{\B,\HH,\SS}(\pi^*\alpha) \\
&= &
\chi^{\B',\HH',\SS'}({\pi'}^*\alpha) \\
&= &
\int_{\B'/T}{\pi'}^*\alpha\wedge e_T(\HH') \\
&= &
\int_{\B'/T}{\pi'}^*\alpha\wedge{\pi'}^*e_T(F\oplus E') \\
&= &
\int_{X/T_0}\Res_\infty\left(
\frac{\alpha\wedge e_T(F\oplus E')}{e_T(E\oplus E')}
\right) \\
&= &
\int_{X/T_0}\Res_\infty
\left(\frac{\alpha}{e_T(\IND(D))}\right).
\end{eqnarray*}
Here the second equation uses the {\it (Naturality)}
axiom for the Euler class, the third equation uses
the {\it (Thom class)} axiom (see~\cite{CMS}),
the fourth equation uses the fact that $\HH'$ is the pullback
of $F\oplus E'$, and the fifth equation uses the hypothesis
that the result holds when $E$ is a trivial bundle
and $F=0$.

\medskip
\noindent{\bf Step~2.}
{\it Suppose $E=L=X\times\C$ is a trivial line bundle
and denote by $B\subset L$ the unit circle bundle.
Then for every $\alpha=\sum_{j\le m/2}\alpha_jz^j\in\Om^m_{T,T_1}(X)$,
\begin{equation}\label{eq:local1}
     \int_{B/T}\pi^*\alpha_0
     = \int_{X/T_0}\Res_\infty\left(\frac{\alpha}{e_T(L)}\right).
\end{equation}
}

\medskip
\noindent
Let $\rho:T\to S^1$ denote the covariant action
of $T$ on the fibres of $L$ and suppose that $T_1$ 
acts on the fibers with weight $\ell$. Then
$\dot\rho(e_1)=-2\pi i\ell$ and, by Lemma~\ref{le:euler-conn},
$$
     e_T(L,\xi+ze_1) = \ell z+\frac{i\dot\rho(\xi)}{2\pi}.
$$
Each form $\alpha_j\in\Om^{m-2j}_T(X)$ is equivariantly closed 
and hence represents a $T$-equivariant cohomology class on $X$. 
Now
\begin{eqnarray*}
    \frac{\alpha(\xi,z)}{e_T(\xi+ze_1)}
&= &
    \left(\sum_{j\le m/2}\frac{\alpha_j(\xi)}{\ell}z^{j-1}\right)
    \frac{1}{1+\frac{i\dot\rho(\xi)}{2\pi\ell z}}  \\
&= &
    \left(\sum_{j\le m/2}\frac{\alpha_j(\xi)}{\ell}z^{j-1}\right)
    \sum_{k\ge 0}\left(-\frac{i\dot\rho(\xi)}{2\pi\ell z}\right)^k.
\end{eqnarray*}
The residue at infinity is minus the coefficient of $z^{-1}$
in this power series.  Thus
\begin{equation}\label{eq:res}
    \Res_\infty\left(\frac{\alpha}{e_T(L)}\right)
    = -\frac{1}{\ell}\sum_{j\ge 0}
      \alpha_j(\xi)\left(-\frac{i\dot\rho(\xi)}{2\pi\ell}\right)^j
\end{equation}
By~(\ref{eq:alpha-1}), the right hand side 
is invariant under the shift $\xi\mapsto\xi+te_1$ and 
hence descends to a $T_0$-equivariant differential 
form on $X$. 

To compute the integral of $\pi^*\alpha_0$ over $B/T$,
we denote the elements of $B=X\times S^1$ by $(x,u)$,
where $x\in X$ and $u\in S^1$.  Then $u^{-1}du$
is the standard $\U(1)$-connection form on $B$.
Define $A_1\in\Om^1(B,\tt_1)$ by
$$
     A_1 := -\frac{iu^{-1}du}{2\pi\ell}e_1.
$$
This form is $T$-invariant and satisfies
$A_1(0,-\dot\rho(e_1)z)=(i\dot\rho(e_1)/2\pi\ell)e_1=e_1$.
Hence it is a $T$-invariant $T_1$-connection as
in~\cite{CMS}, where we regard the action by
$\rho^{-1}$ as the contravariant action on $B$.
The infinitesimal covariant
action of $\xi\in\tt$ on $B$ is given by
$X_\xi(x,u)=(0,\dot\rho(\xi)u)$. Hence the 
$T$-equivariant curvature of $A_1$ is the 
$2$-form $F_{A_1,T}\in\Om^2_T(B,\tt)$ given by
$$
F_{A_1,T}(\xi)
:= F_{A_1} + \xi + A_1(X_\xi)
 = \xi - \frac{i\dot\rho(\xi)}{2\pi\ell}e_1.
$$
(see~\cite[Section~3]{CMS}). 
Replacing $\xi$ by $F_{A_1,T}(\xi)$ in the 
equivariant differential form $\pi^*\alpha_0\in\Om^m_T(B)$
we obtain the $T_1$-basic $T$-equivariant differential
form
$$
(\pi^*\alpha_0)_{A_1}
= \pi^*\alpha_0\left(\xi - \frac{i\dot\rho(\xi)}{2\pi\ell}e_1\right)
= \sum_{j\ge0}\pi^*\alpha_j(\xi)
  \left(-\frac{i\dot\rho(\xi)}{2\pi\ell}\right)^j
$$
on $B=X\times S^1$. The projection $\pi:B\to X$ induces
a $T_0$-equivariant diffeomorphism from $B/T_1$ to $X$,
however, each point in $B$ has an isotropy subgroup of order $\ell$
under the action of $T_1$. Moreover, the diffeomorphism
is orientation preserving if and only if $\ell$ is
negative. (If $\xi_1,\dots,\xi_m$ is a positive basis 
of $T_xX/T_0$ and $u\in S^1$ then, according to 
Remark~\ref{rmk:orient}, the basis $(\xi_1,0),\dots,(\xi_m,0)$
of $T_{(x,u)}B/T$ is positive if and only if the vectors
$(\xi_1,0),\dots,(\xi_m,0),(0,u),(0,\dot\rho(e_1)u)$
form a positive basis of $T_xX/T_0\times\C$. Since 
$\dot\rho(e_1)=-2\pi i\ell$, this is the case if and only
if $\ell$ is negative.) Hence
$$
    \int_{B/T}\pi^*\alpha_0
    = \int_{B/T}(\pi^*\alpha_0)_{A_1}
    = -\frac{1}{\ell}\sum_{j\ge0}\int_{X/T_0}
       \alpha_j(\xi_0)\left(-\frac{i\dot\rho(\xi_0)}{2\pi\ell}\right)^j,
$$
and so 
the
assertion of Step~2 follows from~(\ref{eq:res}). 

\smallbreak

\medskip
\noindent{\bf Step~3.}
{\it We prove the theorem.}

\medskip
\noindent
By Step~1, we may assume without loss of generality
that $F=0$ and $E=X\times V$ for some unitary
$T$-representation~$V$.  Since $T$ is a torus,
we may assume that $V=\C^n$ and that $T$ acts diagonally
by homomorphisms $\rho_\nu:T\to S^1$ for $\nu=1,\dots,n$.
Denote by $L_\nu:=X\times\C$ the $T$-equivariant bundle
where $T$ acts by $\rho_\nu$ on the fibre.
Consider the $T$-moduli problem $(B,H,S)$ given by
$$
     B:=X\times S^{2n-1},\qquad
     H:=\pi^*L_1\oplus\cdots\oplus\pi^*L_{n-1},
$$
where $\pi:B\to X$ denotes the projection, and
$$
     S(x,z_1,\dots,z_n) := (z_1,\dots,z_{n-1}).
$$
Then the $T$-equivariant Euler class of $H$
is the pullback under $\pi$ of the Euler class
of $L_1\oplus\cdots\oplus L_{n-1}$, i.e.
$$
     e_T(H) = \pi^*e_T(L_1)\wedge\cdots\wedge\pi^*e_T(L_{n-1}).
$$
Let $\alpha\in\Om^{m+2n-2}_T(X)$ and define 
$\beta\in\Om^m_{T,T_1}(X)$ by 
$$
     \beta:=\sum_{j\le m/2}\beta_jz^j
     :=\frac{\alpha}{e_T(L_1)\cdots e_T(L_{n-1})},\qquad
     \beta_j\in\Om^{m-2j}_T(X).
$$
Since $T$ acts on $B$ with finite isotropy, we can represent 
the equivariant cohomology class $e_T(H)\in H^{2n-2}_T(B)$
by a $T$-invariant and horizontal differential form
$\tau_H\in\Om^{2n-2}(B)$ (see~\cite[Theorem~3.8]{CMS}).  
With such a representative the identity 
$\pi^*\alpha=\pi^*\beta\wedge e_T(H)$ 
in $\Om^*_{T,T_1}(X)$ takes the form 
$
     \pi^*\alpha = \pi^*\beta_0\wedge\tau_H.
$
Now $S$ is transverse to the zero section and
$S^{-1}(0)$ is the unit sphere bundle $B_n:=B\cap L_n$ 
in $L_n$. Hence it follows from the {\it (Transversality)}
axiom for the Euler class in~\cite{CMS} that
$$
\int_{B/T}\pi^*\alpha
= \int_{B/T}\pi^*\beta_0\wedge \tau_H  
= \int_{B_n/T}\pi^*\beta_0 
= \int_{X/T_0}\Res_\infty\left(\frac{\beta}{e_T(L_n)}\right).
$$
The last equation follows from Step~2.
Since $\beta/e_T(L_n)=\alpha/e_T(E)$, 
this proves the theorem.
\end{proof}


\subsection{The index formula}\label{sec:index}

We return to the setting of Section~\ref{sec:wall}.
Recall that $P\to\Sigma$ is a principal $T$-bundle
and $L_\nu=P\times_{\rho_\nu}\C\to\Sigma$ for $\nu=1,\dots,n$,
Given an index set $I\subset\{1,\dots,n\}$ as in 
Section~\ref{sec:wall} we consider the principal 
$\Gg_0$-bundle
$$
     \Pp_I := \Aa(P)\times\bigoplus_{\nu\in I}\Om^0(\Sigma,L_\nu)
     \to\Bb_I := \Pp_I/\Gg_0
$$
where the based gauge group $\Gg_0$ acts by
$g^*(A,u_\nu)=(A+g^{-1}dg,\rho_\nu(g)^{-1}u_\nu)$.
It also acts on $L_{\nu_0}$ by
$
    g^*[p,\zeta]
    := [pg(z)^{-1},\zeta]
     = [p,\rho_{\nu_0}(g(z))^{-1}\zeta],
$
where $z:=\pi(p)\in\Sigma$.
For $\nu_0\notin I$ we consider the universal line bundle
$$
    \LL^{\nu_0}:=\frac{\Pp_I\times L_{\nu_0}}{\Gg_0}
    \to\Bb_I\times\Sigma.
$$
The torus $T$ acts on $\LL^{\nu_0}$ by
$
    h^*[A,u_\nu,p,\zeta]
    := [A,\rho_\nu(h)^{-1}u_\nu,p,\rho_{\nu_0}(h)^{-1}\zeta].
$

For ${\bf x}=[A,u]\in\Bb_I$ let us denote by
$\LL^\nu_{\bf x}$ the restriction of
$\LL^\nu$ to $\{{\bf x}\}\times\Sigma$.
This restriction is equipped with a natural
connection (induced by $A$) and hence
with a Cauchy-Riemann operator
$$
     \bar\p^\nu_{\bf x}:
     \Om^0(\Sigma,\LL^\nu_{\bf x})\to
     \Om^{0,1}(\Sigma,\LL^\nu_{\bf x}).
$$
Next consider the universal vector bundle
$$
     \E := \bigoplus_{\nu\notin I}\LL^\nu
$$
and its restrictions $\E_{\bf x}$ to $\{{\bf x}\}\times\Sigma$.
The Cauchy-Riemann operators
$$
     \bar\p_{\bf x}:\Om^0(\Sigma,\E_{\bf x})\to
     \Om^{0,1}(\Sigma,\E_{\bf x})
$$
form a family of Fredholm operators over $\Bb_I$
between appropriate Hilbert space completions
$\Ee_{\bf x}$ of $\Om^0(\Sigma,\E_{\bf x})$
and $\Ff_{\bf x}$ of $\Om^{0,1}(\Sigma,\E_{\bf x})$.
These operators are complex linear and equivariant
with respect to the action of $T$.

As in Section~\ref{sec:wall} we denote by $T_1$
the identity component of the isotropy subgroup
of the subspace
$\C^I:=\left\{x\in\C^n\,|\,x_\nu=0\mbox{ for }\nu\notin I\right\}$
and assume that $T_1$ is a circle.
This circle acts trivially on the base $\Bb_I$
and with finite isotropy outside of the zero sections
of $\Ee$ and $\Ff$. The quotient group $T_0:=T/T_1$ acts
with finite isotropy on the moduli space $\Mm_0\subset\Bb_I$
of solutions of equation~(\ref{eq:v0}).
Hence we are in the situation of Theorem~\ref{thm:local}.
The relevant dimensions are
\begin{eqnarray*}
\dim\Mm - \dim T
&= &
(n-\dim T)(2-2g) + 2\sum_{\nu=1}^nd_\nu =:2m, \\
\dim\Mm_0 - \dim T_0
&= &
(|I|-\dim T_0)(2-2g) + 2\sum_{\nu\in I}d_\nu + 2g, \\
\INDEX(\bar\p)
&= &
\sum_{\nu\notin I}(d_\nu+1-g),
\end{eqnarray*}
where $d_\nu$ denotes the first Chern number of the bundle
$L_\nu\to\Sigma$. Note that the kernel manifold of
$\bar\p$ is precisely the space $\Ss_0$ of solutions of
equations~(\ref{eq:vs}). Hence, 
by Theorem~\ref{thm:local} and Remark~\ref{rmk:local}, 
we have 
\begin{eqnarray}\label{eq:local3}
\int_{\Ss_0/T}\pi^*\alpha
&= &
\int_{\Mm_0/T_0}\Res_\infty\left(
\frac{\alpha}{e_T(\IND(\bar\p))}
\right) 
\nonumber \\
&= &
\int_{\Mm_0/T_0}\Res_\infty\left(
\frac{\alpha}{\prod_{\nu\notin I}e_T(\IND(\bar\p^\nu))}
\right).
\end{eqnarray}

\begin{lemma}\label{le:index}
Denote by $\pi_0^*:H_{T,T_1}^*(\Aa/\Gg_0)\to H^*_{T,T_1}(\Bb_I)$
the homomorphism induced by the projection
$\pi_0:\Bb_I\to\Aa/\Gg_0$. Then, for every $\nu\notin I$, 
$$
     e_T(\IND(\bar\p^\nu))
     = \pi_0^*\inner{\w_\nu}{\xi}^{d_\nu+1-g}
       \exp\left(-\frac{\Om_\nu}{\inner{\w_\nu}{\xi}}\right)
     \in H^*_{T,T_1}(\Bb_I).
$$
\end{lemma}

\begin{proof}[Proof of Theorem~\ref{thm:WCN}]
By Proposition~\ref{prop:wcn} and~(\ref{eq:local3})
we have
\begin{eqnarray*}
\Phi^{\rho,\tau_0+\eps\tau_1}_{\lambda,g}(\alpha)
- \Phi^{\rho,\tau_0-\eps\tau_1}_{\lambda,g}(\alpha)
&= &
-\int_{\Ss_0/T}\pi^*\alpha \\
&= &
-\int_{\Mm_0/T_0}\Res_\infty\left(\frac{\alpha}
{\prod_{\nu\notin I}e_T(\IND(\bar\p^\nu))}\right)  \\
&= &
\int_{\Mm_0/T_0}\pi_0^*\alpha_0,
\end{eqnarray*}
where
$$
\alpha_0(\xi)
= 
\frac{1}{2\pi i}\oint
\frac{\alpha(\xi+ze_1)}
{\prod_{\nu\notin I}\inner{\w_\nu}{\xi+ze_1}^{d_\nu+1-g}}
\exp\left(\sum_{\nu\notin I}\frac{\Om_\nu}{\inner{\w_\nu}{\xi+ze_1}}
\right)\,dz
$$
(by Lemma~\ref{le:index}). 
This proves the theorem. 
\end{proof}

\begin{proof}[Proof of Lemma~\ref{le:index}]
Note first that $\LL^\nu$ is the pullback under the projection 
$\pi_0\times\id:\Bb_I\times\Sigma\to\Aa/\Gg_0\times\Sigma$ 
of the bundle
$$
\Ll^\nu := \frac{\Aa\times L_\nu}{\Gg_0}
\to \Aa/\Gg_0\times\Sigma,
$$
and $\IND(\bar\p^\nu)\in K_T(\Bb_I)$ is the pullback under $\pi_0$
of the index bundle of the Cauchy-Riemann operators on $\Ll^\nu$.
The torus $T$ acts trivially on $\Aa/\Gg_0\times\Sigma$
and by $\rho_\nu$ on the fibres of $\Ll^\nu$. Hence by
Example~\ref{ex:euler-diag},
\begin{equation}\label{eq:eT}
e_T(\IND(\bar\p^\nu,\Ll^\nu))
= \sum_{j\ge0}
  \inner{\w_\nu}{\xi}^{d_\nu+1-g-j}
  c_j(\IND(\bar\p^\nu,\Ll^\nu)).
\end{equation}
Hence it remains to compute the {\em ordinary} Chern classes 
of the $K$-theory class $\IND(\bar\p^\nu,\Ll^\nu)$.
The Atiyah--Singer index theorem for families asserts that
$$
\ch(\IND(\bar\p^\nu,\Ll^\nu))
= \int_\Sigma \td(T\Sigma)\ch(\Ll^\nu)
\in H^*(\Aa/\Gg).
$$
(See~\cite[Theorem~4.3]{AS1} and~\cite[Theorem~5.1]{AS2}.)
Here $\ch$ and $\td$ denote the Chern character and
the Todd class, respectively.
The Todd class of $T\Sigma$ is given by
$$
\td(T\Sigma) = 1 + (1-g)\sigma,
$$
where $\sigma\in H^2(\Sigma;\Z)$ denotes the positive generator. Thus
our task at hand is to compute the Chern character of the line
bundle $\Ll^\nu$.  By Lemma~\ref{lem:c1} below, the first 
Chern class of $\Ll^\nu$ is given by
$$
c_1(\Ll^\nu) = d_\nu\sigma 
- \sum_{i=1}^{2g}\sum_{j=1}^k \w_{\nu j}\alpha_i\wedge\tau_{ij}.
$$
{From} this we can compute $\ch(\Ll^\nu)$. Note that
$\alpha_i\wedge\alpha_{i'}=\pm\sigma$ whenever $i'=i\pm g$
and is equal to zero otherwise. Hence
\begin{eqnarray*}
     \frac12\left(\sum_{i=1}^{2g}\sum_{j=1}^k
     \w_{\nu j}\alpha_i\wedge\tau_{ij}\right)^2
&= &
     -\sigma\wedge\sum_{i=1}^g\sum_{j,j'=1}^k
     \w_{\nu j}\w_{\nu j'}\tau_{ij}\wedge\tau_{i+g,j'} \\
&=: &
     -\sigma\wedge\Om_\nu,
\end{eqnarray*}
and all higher powers vanish. It follows that
$$
\ch(\Ll^\nu) 
= 1 + d_\nu\sigma 
  - \sum_{i=1}^{2g}\sum_{j=1}^k
    \w_{\nu j}\alpha_i\wedge\tau_{ij} 
  - \sigma\wedge\Om_\nu.
$$
Applying the index theorem for families
we obtain
\begin{eqnarray*}
\ch(\IND(\bar\p^\nu,\Ll^\nu))
&= &
\int_\Sigma \td(T\Sigma)\ch(\Ll^\nu) \\
&= &
\int_\Sigma \left(
(d_\nu+1-g)\sigma
- \sum_{i=1}^{2g}\sum_{j=1}^k
\w_{\nu j}\alpha_i\wedge\tau_{ij}
- \sigma\wedge\Om_\nu
\right) \\
&= &
d_\nu+1-g - \Om_\nu.
\end{eqnarray*}
The last formula implies, by a standard algebraic argument,
that
$$
c_1(\IND(\bar\p^\nu,\Ll^\nu)) = -\Om_\nu,\qquad
c_j(\IND(\bar\p^\nu,\Ll^\nu))
= \frac{1}{j!}c_1(\IND(\bar\p^\nu,\Ll^\nu))^j.
$$
Hence, by~(\ref{eq:eT}),
\begin{eqnarray*}
     e_T(\IND(\bar\p^\nu,\Ll^\nu))
&= &
     \sum_{j\ge0}
     \inner{\w_\nu}{\xi}^{d_\nu+1-g-j}
     \frac{(-\Om_\nu)^j}{j!} \\
&= &
     \inner{\w_\nu}{\xi}^{d_\nu+1-g}
     \exp\left(-\frac{\Om_\nu}{\inner{\w_\nu}{\xi}}\right).
\end{eqnarray*}
Since $e_T(\IND(\bar\p^\nu))=\pi_0^*e_T(\IND(\bar\p^\nu,\Ll^\nu))$
the result follows.
\end{proof}

\begin{lemma}\label{lem:c1}
The first Chern class of $\Ll^\nu$ is given by
$$
c_1(\Ll^\nu) = d_\nu\sigma -
\sum_{i=1}^{2g}\sum_{j=1}^k \w_{\nu j}\alpha_i\wedge\tau_{ij}.
$$
\end{lemma}

\begin{proof}
Fix a reference connection $A_0\in\Aa(P)$ and denote by 
$\Aa_0\subset\Aa(P)$ the set of connections $A\in\Aa(P)$ 
that satisfy $F_A=\Vol(\Sigma)^{-1}\lambda(P)\dvol_\Sigma$ 
and $d^*(A-A_0)=0$. The restricted gauge group 
$\Gg_{00}\subset\Gg_0(P)$ consists of all gauge transformations 
$g:\Sigma\to T$ that satisfy $d^*(g^{-1}dg)=0$ and 
$g(z_0)=\one$.
Identify the quotient $\Aa_0/\Gg_{00}$ with the 
standard torus $\T^{2gk}$ via the map
$$
      \R^{2gk}\to\Aa_0:t\mapsto A_t
      := A_0 + \sum_{i=1}^{2g}\sum_{j=1}^k t_{ij}\alpha_i e_j.
$$
For $i$ and $j$ let $g_{ij}\in\Gg_{00}$ be the unique 
harmonic gauge transformation that satisfies 
$g_{ij}^{-1}dg_{ij}=\alpha_ie_j$ 
(and $g_{ij}(z_0)=\one$). 
Then the restriction of $\Ll^\nu$ to the submanifold 
$\Aa_0/\Gg_{00}\times\Sigma$ can be identified with the 
quotient $\R^{2gk}\times L_\nu/\Z^{2gk}$, where 
$m=\{m_{ij}\}\in\Z^{2gk}$ acts by 
$$
     m^*(t,z,v):=(t+m,z,\prod_{ij}\rho_\nu(g_{ij}(z))^{-m_{ij}}v).
$$
A section is a map
$\R^{2gk}\times\Sigma\to L_\nu:(t,z)\mapsto s(t,z)=s_t(z)\in L_{\nu z}$ 
that satisfies
$$
     s(t+m,z) = \prod_{ij}\rho_\nu(g_{ij}(z))^{-m_{ij}}s(t,z),\qquad
     m\in\Z^{2gk}. 
$$
A connection is given by the formula 
$$   
     d^\nabla s := d_{A_t}s_t 
     + \sum_{i=1}^{2g}\sum_{j=1}^k \frac{\p s_t}{\p t_{ij}}dt_{ij}.
$$
An easy computation shows that the curvature of this connection 
is the $2$-form $F^\nabla\in\Om^2(\T^{2gk}\times\Sigma,\sqrt{-1}\R)$
given by 
\begin{eqnarray*}
     F^\nabla 
&= &
     \dot\rho_\nu(F_{A_t})
     - \sum_{i=1}^{2g}\sum_{j=1}^k\dot\rho_\nu
     \left(\frac{\p A_t}{\p t_{ij}}\right)\wedge dt_{ij}  \\
&= &
     -2\pi\sqrt{-1}\inner{\w_\nu}{\lambda(P)}
     \frac{\dvol_\Sigma}{\Vol(\Sigma)}
     + \sum_{i=1}^{2g}\sum_{j=1}^k
       2\pi\sqrt{-1}\inner{\w_\nu}{e_j}\alpha_i\wedge dt_{ij}.
\end{eqnarray*}
Since the first Chern class of of $\Ll^\nu$ is represented
by the $2$-form $\sqrt{-1}F^\nabla/2\pi$, the result follows. 
\end{proof}


\section{Computation of the genus zero invariants}\label{sec:compute}

Let $\Ii$ denote the set of partitions 
$$
      I_1\cup\cdots\cup I_k = \{1,\dots,n\} 
$$
that satisfy
the (Dimension) condition in the introduction, i.e.~for every 
$j\in\{1,\dots,k\}$, the subspace 
$$
      E_j:=\SPAN\{\w_\nu\,|\,\nu\in I_1\cup\cdots\cup I_j\}
      \subset\tt^*
$$
has dimension $j$ and $\w_\nu\notin E_j$ for every 
$\nu\in I_{j+1}\cup\cdots\cup I_k$.  
It is now convenient to fix an orientation of $\tt$. 
For every $I=(I_1,\dots,I_k)\in\Ii$ we introduce 
the iterated residue
$
     \phi^I_\lambda=\phi^{\rho,I}_\lambda:
     S^*(\tt^*)\to\R
$
by
\begin{equation}\label{eq:residue}
        \phi_\lambda^I(\alpha)
        := \frac{1}{(2\pi i)^k}\oint\dots\oint \frac{\alpha(\sum z_je_j)}
        {\prod_{\nu=1}^n\inner{\w_\nu}{\sum z_je_j}^{d_\nu+1}}
        dz_k\cdots dz_1,
\end{equation}
where $d_\nu:=\inner{\w_\nu}{\lambda}$.
Here the lattice vectors $e_1,\dots,e_k\in\tt$ form an oriented 
integral basis of $\Lambda$ such that
the vectors $e_j,\dots,e_k$ are orthogonal to the span of 
the vectors $\w_\nu$ for $\nu\in I_1\cup\cdots\cup I_{j-1}$ 
and $2\leq j\leq k$.  These requirements determine the $e_j$ 
up to a change $e_j' = \pm e_j+\sum_{i>j}a_{ij}e_i$. 
The corresponding coordinates $\xi=\sum z_je_j=\sum z_j'e_j'$ 
change by $z_i=\pm z_i'+\sum_{j<i}a_{ij}z_j'$. 
Since the $e_i$ and the $e_i'$ form oriented bases 
there is an even number of minus signs. It follows from these
observations and Remark~\ref{rmk:residue} that the 
integral~(\ref{eq:residue}) is independent of the choice
of the $e_i$.

\begin{lemma}\label{le:sum}
For every regular value $\tau$ of $\mu$ there exists
a collection of integers $\{m_I\}_{I\in\Ii}$ such that
$$
     \Phi^{\rho,\tau}_\lambda(\alpha)
     =\sum_{I\in\Ii}m_I\phi^I_\lambda(\alpha)
$$
for every $\lambda\in\Lambda$ and every $\alpha\in S^*(\tt^*)$.
Moreover, given a path $\gamma$ connecting $\tau$ to $\tt^*\setminus
\im\,\mu$, the coefficients $m_I$ can be chosen to satisfy the
following condition: If $I=(I_1,\dots,I_k)\in\Ii$ is such that 
$\gamma$ does not intersect the hyperplane spanned 
by the $\w_\nu$ for $\nu\in I_1\cup\cdots\cup  I_{k-1}$, then $m_I=0$.
\end{lemma}

\begin{proof}
Theorem~\ref{thm:wall0} and induction over the dimension of $T$. 
\end{proof}

\begin{remark}\label{rmk:hyperplane}\rm
Fix an oriented basis $e_1,\dots,e_k$ of $\tt$,
let $H_\nu\subset\C^k$ be the hyperplane 
$\sum_{j=1}^k\inner{\w_\nu}{e_j}z_\nu=0$,
and denote
$
     H := \bigcup_{\nu=1}^n H_\nu.
$
Then the right hand side of~(\ref{eq:residue})
can be interpreted as the integral of the 
$k$-form 
$$
     \om_{\lambda,\alpha}
     := \frac{\alpha(\sum z_je_j)}
     {(2\pi i)^k\prod_{\nu=1}^n
     \inner{\w_\nu}{\sum z_je_j}^{\inner{\w_\nu}{\lambda}+1}}
     dz_1\wedge\cdots\wedge dz_k
     \in\Om^k(\C^k\setminus H)
$$
over a suitable homology class $\sigma_I\in H_k(\C^k\setminus H)$.
Hence, by Lemma~\ref{le:sum},
there is a locally constant 
map $\tau\mapsto\sigma(\tau):=\sum_{I\in\Ii}m_I(\tau)\sigma_I$
which assigns to every regular value of $\mu$ a homology
class $\sigma(\tau)\in H_k(\C^k\setminus H)$ such that
the invariant $\Phi^{\rho,\tau}_\lambda(\alpha)$
is equal to the integral of $\om_{\lambda,\alpha}$ 
over $\sigma(\tau)$ (for all $\lambda$ and $\alpha$).
It is an interesting problem to study the map 
$\tau\mapsto\sigma(\tau)$ in more detail.
\end{remark}

\begin{lemma}\label{le:-1}
For $\nu=1,\dots,n$ let $\ell_\nu$ be integers and $a_\nu,b_\nu$ 
be real numbers such that $a_\nu\neq 0$.
If $\sum_{\nu=1}^n\ell_\nu=-1$ then
$$
     \frac{1}{2\pi i}\oint
     \prod_{\nu=1}^n(a_\nu z+b_\nu)^{\ell_\nu}dz
     = \prod_{\nu=1}^na_\nu^{\ell_\nu}.
$$
If $\sum_{\nu=1}^n\ell_\nu<-1$ then the integral is zero.
\end{lemma}

\begin{proof}
In the variable $w:=1/z$ the integrand reads
$$
        \prod_\nu\left(\frac{a_\nu}{w} +
        b_\nu\right)^{\ell_\nu}\left(\frac{-dw}{w^2}\right)  
        = -\frac{\prod_\nu(a_\nu+b_\nu w)^{\ell_\nu}dw}
        {w^{2+\sum\ell_\nu}}. 
$$
Since $a_\nu\neq 0$ the numerator is holomorphic near the 
origin. Hence the residue is zero whenever 
$2+\sum_\nu\ell_\nu\le 0$ and is equal to 
$\prod_\nu a_\nu^{\ell_\nu}$ whenever $2+\sum_\nu\ell_\nu=1$.
\end{proof}

\begin{lemma}\label{le:well}
Let $I\in\Ii$ and $e_1,\dots,e_k$ be a positive basis of 
$\Lambda$ such that $e_j,\dots,e_k$ are orthogonal to $\w_\nu$ 
for $\nu\in I_1\cup\cdots\cup I_{j-1}$ (as in~(\ref{eq:residue})). 
Let $\lambda\in\Lambda$ and $\ell$ be an $n$-tuple of 
nonnegative integers 
such that 
\begin{equation}\label{eq:dim}
     |\ell| = n-k+\sum_{\nu=1}^nd_\nu,\qquad
     d_\nu:=\inner{\w_\nu}{\lambda}.
\end{equation}
If $I\in\Ii_\lambda(\ell)$ then 
$$
        \phi_\lambda^I(\w^\ell) =
        \prod_{j=1}^k\prod_{\nu\in I_j}
        \inner{\w_\nu}{e_j}^{\ell_\nu-d_\nu-1}.
$$ 
Otherwise $\phi_\lambda^I(\w^\ell)=0$. 
\end{lemma}

\begin{proof}
The condition $I\in\Ii_\lambda(\ell)$ asserts that 
$$
     \sum_{\nu\in I_j}(\ell_\nu-d_\nu-1)=-1
$$
for $j=1,\dots,k$.  Consider the integral over $z_k$. 
The coefficient $\inner{\w_\nu}{e_k}$ of $z_k$ in the linear 
map $(z_1,\dots,z_k)\mapsto\inner{\w_\nu}{\sum_jz_je_j}$ 
is nonzero iff $\nu\in I_k$. So, by Lemma~\ref{le:-1}, 
\begin{eqnarray*}
     \phi^I_\lambda(\w^\ell)
&= &
     \left(\prod_{\nu\in I_k}
     \inner{\w_\nu}{e_k}^{\ell_\nu-d_\nu-1}\right)
     \frac{1}{(2\pi i)^{k-1}}  \\
&&
        \oint\cdots\oint\prod_{j=1}^{k-1}
        \prod_{\nu\in I_j}
        \inner{\w_\nu}{\sum_{j=1}^{k-1}z_je_j}^{\ell_\nu-d_\nu-1}
        dz_{k-1}\dots dz_1 
\end{eqnarray*}
whenever $\sum_{\nu\in I_k}(\ell_\nu-d_\nu-1)=-1$. 
If $\sum_{\nu\in I_k}(\ell_\nu-d_\nu-1)<-1$  
then the integral over $z_k$ is zero. 
Hence it follows by induction that 
$\phi^I_\lambda(\w^\ell)$ has the required form
whenever $I\in\Ii_\lambda(\ell)$. 
If $I\notin\Ii_\lambda(\ell)$ then it follows
from~(\ref{eq:dim}) that $\sum_{\nu\in I_j}(\ell_\nu-d_\nu-1)<-1$
for some $j$ and hence $\phi^I_\lambda(\w^\ell)=0$.
\end{proof}

\begin{proof}[Proof of Theorem~\ref{thm:invariants}.]
Assertion~(i) follows from Lemmata~\ref{le:sum}
and~\ref{le:well}.

We prove~(ii). Let $\lambda\in\Lambda$,
$d_\nu:=\inner{\w_\nu}{\lambda}$, $\ell$ be an $n$-tuple of 
nonnegative integers, $J\subset\{1,\dots,n\}$
be an index set such that $\{\w_\nu\,|\,\nu\in J\}$
is a basis of $\tt^*$, and assume that $\ell_\nu=d_\nu$
for $\nu\in J$ and $\ell_\nu=d_{\nu+1}$ for $\nu\notin J$. 
Then a partition $I$ belongs to $\Ii_\lambda(\ell)$ 
if and only if $I\in\Ii$ and $I_j\cap J$ consists of a single 
element for each $j$. This follows from the fact that 
$\sum_{\nu\in I_j}(\ell_\nu-d_\nu-1)=-1$ and that each 
summand $\ell_\nu-d_\nu-1$ equals $0$ or $-1$.

Assume $\tau\notin C(J)
=\{\sum_{\nu\in J}\eta_\nu \w_\nu\;|\;\eta_\nu\geq 0\}$.  
We must prove that $\Phi^{\rho,\tau}_\lambda(\w^\ell)=0$.
To see this, we examine the set $\Ii_\lambda(\ell)$.
Since the set $\{\w_\nu\,|\,\nu\in J\}$ is linearly independent 
it follows that, for each ordering $J=\{\nu_1,\dots,\nu_k\}$, 
there exists a unique partition $I\in\Ii_\lambda(\ell)$ 
such that $\nu_j\in I_j$ for all $j$, and conversely,
each partition $I\in\Ii_\lambda(\ell)$ determines an
ordering of $J$.  Moreover, for every such partition the 
hyperplane $W_I:=\SPAN\{\w_\nu\,|\,\nu\notin I_k\}$
agrees with the hyperplane 
$W_{\nu_k}:=\SPAN\{\w_\nu\,|\,\nu\in J\setminus\{\nu_k\}\}$.
Hence the hyperplanes $\W_I$ for $I\in\Ii_\lambda(\ell)$
are precisely the supporting hyperplanes of $C(J)$.
Since $\tau\notin C(J)$, there exists a straight line $\gamma$ 
connecting $\tau$ to $\tt^*\setminus\im\,\mu$ which misses the 
supporting hyperplanes (this is true because $C(J)$ is a cone over
a simplex).  Hence the coefficients $m_I$ in Lemma~\ref{le:sum} can be 
chosen such that $m_I=0$ for every $I\in\Ii_\lambda(\ell)$. 
This implies that
$$
     \Phi^{\rho,\tau}_\lambda(\w^\ell) 
     = \sum_{I\in\Ii_\lambda(\ell)}m_I\phi_\lambda^I(\w^\ell) 
     = 0.
$$
The same argument shows that the invariant
$\Phi^{\rho,\tau}_\lambda(\w^\ell)$ for $\tau\in C(J)$
is independent of $\tau$. 

Assume $\tau\in C(J)$.  
Fix an ordering $J=\{\nu_1,\dots,\nu_k\}$
and let $I\in\Ii_\lambda(\ell)$ be the 
unique partition satisfying $\nu_j\in I_j$ for $j=1,\dots,k$.
Choose an integer basis $e_1,\dots,e_k$ of $\tt$
such that $\inner{\w_\nu}{e_j}=0$ for 
$\nu\in I_1\cup\cdots I_{j-1}$ and 
$
     \inner{\w_{\nu_j}}{e_j}>0.
$
Let $\tau_0$ be a positive linear 
combination of $\w_{\nu_1},\dots,\w_{\nu_{k-1}}$.
Since the invariant is independent of the choice
of $\tau\in C(J)$, we may assume $\tau=\tau_0+\eps\tau_1$,
where $\tau_1:=\w_{\nu_k}$.  Since the invariant is zero 
outside of $C(J)$, 
we have 
$$
     \Phi^{\rho,\tau_0-\eps\tau_1}_\lambda(\w^\ell)=0.
$$
Hence, by Theorem~\ref{thm:wall0}, 
$$
     \Phi^{\rho,\tau}_\lambda(\w^\ell)
     = \frac{1}{\inner{\w_{\nu_k}}{e_k}}
     \Phi^{\rho_0,\tau_0}_{\lambda_0}\left(
     \prod_{j=1}^{k-1}\prod_{\nu\in I_j}\w_\nu^{\ell_\nu}
     \right). 
$$
Now assertion~(ii) follows by induction. 

We prove~(iii).  
Assume $|\ell|=n-k+\sum_{\nu=1}^nd_\nu$
(otherwise both invariants are zero).
Since
$$
     \ell_\nu-d_\nu-1
     = \ell_\nu+d_\nu'-(d_\nu+d_\nu')-1,
$$
for every $\nu$ we have 
$$
     \Ii_\lambda(\ell)=\Ii_{\lambda+\lambda'}(\ell+d'),\qquad
     \phi^I_\lambda(\w^\ell) = \phi^I_{\lambda+\lambda'}(\w^{\ell+d'})
$$
for every $I\in\Ii$ (see Lemma~\ref{le:well}). 
Hence~(iii) follows from Lemma~\ref{le:sum}. 

To prove~(iv) and~(v) we introduce the following notation.
For every $\lambda\in\Lambda$ and every $n$-tuple
$\ell$ of nonnegative integers define the number 
$$
     \i_\lambda(\ell) 
     := \sum_{d_\nu\ge0}\max\{\ell_\nu-d_\nu-1,0\}
        + \sum_{d_\nu\le-1}\ell_\nu,\qquad
     d_\nu := \inner{\w_\nu}{\lambda}.
$$
Recall that $J_\ell:=\left\{\nu\,|\,\ell_\nu\le d_\nu\right\}$.
We prove~(iv) and~(v) in five steps.

\smallskip
\noindent{\bf Step~1.}
{\it 
Each $\w^\ell$ can be expressed as a linear combination 
of classes $\w^{\ell'}$ that satisfy $J_{\ell'}\subset J_\ell$
and either $\Ii_\lambda(\ell')=\emptyset$ or 
$
      \ell_\nu'\le\max\{d_\nu+1,0\}
$
for all $\nu$.}

\smallskip
\noindent
We prove Step~1 by induction over $\i_\lambda(\ell)$.
If $\i_\lambda(\ell)=0$ then 
$
      \ell_\nu'\le\max\{d_\nu+1,0\}
$
for all $\nu$.
Assume $\i_\lambda(\ell)>0$ and, by induction, that 
the claim has been established for every $\ell'$
that satisfies $\i_\lambda(\ell')<\i_\lambda(\ell)$.
If $\Ii_\lambda(\ell)=\emptyset$ there is nothing to prove.  
Hence assume $\Ii_\lambda(\ell)\ne\emptyset$.  
Since $\i_\lambda(\ell)>0$ there is a $\nu_0$ 
such that 
$
     \ell_{\nu_0}>0
$
and
$
    \ell_{\nu_0}>d_{\nu_0}+1.
$
Let $I=(I_1,\dots,I_k)\in\Ii_\lambda(\ell)$.
Since $\sum_{\nu\in I_j}(\ell_\nu-d_\nu-1)=-1$
for every $j$ there are indices $\nu_j\in I_j$
for $j=1,\dots,k$ such that $\ell_{\nu_j}\le d_{\nu_j}$.
By the (Dimension) condition, the vectors $\w_{\nu_1},\dots,\w_{\nu_k}$
form a basis of $\tt^*$.  Hence $\w_{\nu_0}$ can be expressed
as a linear combination of the vectors $\w_{\nu_j}$. 
Since $\ell_{\nu_0}>0$ we can replace one of the factors 
$\w_{\nu_0}$ in $\w^\ell$ by this linear combination. 
This expresses $\w^\ell$ as a linear combination of monomials 
of the form $\w^{\ell'}$ with $\i_\lambda(\ell')<\i_\lambda(\ell)$
and $J_{\ell'}\subset J_\ell$.  Hence the assertion for 
$\w^\ell$ follows from the induction hypothesis.  

\smallskip
\noindent{\bf Step~2.}
{\it 
If $\ell_\nu\le d_\nu+1$ for every $\nu$
and $\Ii_\lambda(\ell)\ne\emptyset$
then $\ell$ satisfies~(ii).}

\smallskip
\noindent
Let $I\in\Ii_\lambda(\ell)$.  Then 
the formula $\sum_{\nu\in I_j}(\ell_\nu-d_\nu-1)=-1$
shows that, for each $j$, there is precisely 
one index $\nu_j\in I_j$ such that $\ell_{\nu_j}=d_{\nu_j}$
and $\ell_\nu=d_\nu+1$ for $\nu\in I_j\setminus\{\nu_j\}$.
Since the vectors $\w_{\nu_1},\dots,\w_{\nu_k}$ form a basis
of $\tt^*$, it follows that $\ell$ satisfies~(ii) 
with $J=\{\nu_1,\dots,\nu_k\}$. 

\smallskip
\noindent{\bf Step~3.}
{\it We prove~(iv).}

\smallskip
\noindent
Assume $d_\nu\ge-1$ for every $\nu$. Then, by Step~1,
each $\w^\ell\in S^{m_\lambda}(\tt^*)$ is a linear 
combination of classes $\w^{\ell'}$ that satisfy 
either $\Ii_\lambda(\ell')=\emptyset$ or 
$
      \ell_\nu'\le d_\nu+1
$
for all $\nu$. Hence the assertion follows from Step~2.

\smallskip
\noindent{\bf Step~4.}
{\it (v) holds under the assumption $d_\nu\ge-1$ for all $\nu$.}

\smallskip
\noindent
We argue indirectly and assume that 
$\Phi^{\rho,\tau}_\lambda(\w^\ell)\ne0$. 
Then the linear combination in Step~1 must contain 
a term $\w^{\ell'}$ that satisfies $J_{\ell'}\subset J_\ell$,
and $\Phi^{\rho,\tau}_\lambda(\w^{\ell'})\ne0$.
The latter implies that $\Ii_\lambda(\ell')\ne\emptyset$
and so $\ell'_\nu\le d_\nu+1$ for all $\nu$.
Hence, by Step~2, $\ell'$ satisfies~(ii)
with $J=J_{\ell'}=\{\nu\,|\,\ell_\nu'\le d_\nu\}$.
Since $\Phi^{\rho,\tau}_\lambda(\w^{\ell'})\ne0$,
it follows from~(ii) that $\tau\in C(J_{\ell'})\subset C(J_\ell)$.

\smallskip
\noindent{\bf Step~5.}
{\it We prove~(v).}

\smallskip
\noindent
Suppose $\Phi_\lambda(\w^\ell)\ne 0$. 
Choose $\lambda'$ such that
$d_\nu':=\inner{\w_\nu}{\lambda'}\ge\max\{0,{-1-d_\nu}\}$
for all $\nu$.  Then, by~(iii), we have
$
     \Phi_{\lambda+\lambda'}(\w^{\ell+d'})=\Phi_\lambda(\w^\ell)\ne0.
$ 
Hence, by Step~4, $\tau\in C(J)$, where 
$
    J:=\{\nu\,|\,\ell_\nu+d_\nu'\le\la\w_\nu,\lambda+\lambda'\ra\} =
\{\nu\,|\,\ell_\nu\leq d_\nu\}.
$
This proves the theorem.
\end{proof}


\section{Quantum cohomology}\label{sec:qcoh}

Let $\tau$ be a regular value of $\mu$. Throughout this section we
assume that $T$ acts freely on $\mu^{-1}(\tau)$.
Equivalently, if $J\subset\{1,\dots,n\}$ is an index set 
consisting of $k$ elements such that $\tau\in C(J)$, then 
$\tau$ belongs to the interior of $C(J)$ and the 
determinant of the tuple $\{\w_\nu\,|\,\nu\in J\}$ 
is equal to plus or minus one.  Under this assumption the 
symplectic quotient 
$$
      \bar M:=M\dslash T(\tau)=\mu^{-1}(\tau)/T
$$
is a K\"ahler manifold. We denote 
by $H^*(\bar M)$, respectively $H_*(\bar M)$, 
the quotients of the integral (co)homology groups 
by their torsion subgroups. By Kirwan's theorem, 
the homomorphism $H^*(BT)\to H^*(\bar M)$ is surjective 
and the homomorphism $H_*(\bar M)\to H_*(BT)$ is injective.
For every $\nu$ denote by $\bar\w_\nu\in H^2(\bar M)$ the image of 
the cohomology class $\w_\nu\in\Lambda^*\cong H^2(BT)$
under the Kirwan homomorphism $H^2(BT)\to H^2(\bar M)$.
By Theorem~\ref{thm:coh}, the cohomology class 
$\bar\w_\nu$ vanishes whenever $\tau\notin 
C(\{1,\dots,n\}\setminus\{\nu\})$. The remaining 
classes $\bar\w_\nu$ generate $H^2(\bar M)$.
Hence the image of the homomorphism 
$H_2(\bar M)\to H_2(BT)\cong\Lambda$ is the 
subgroup
$$
     \Lambda(\tau) := \left\{\lambda\in\Lambda\,|\,
     \tau\notin C(\{1,\dots,n\}\setminus\{\nu\})\IMP
     \inner{\w_\nu}{\lambda}=0\right\}.
$$
Recall the definition of the inverse isomorphism 
$\Lambda(\tau)\to H_2(\bar M):\lambda\mapsto\bar\lambda$ and the 
effective cone $\Lambda_\eff(\tau)\subset\Lambda(\tau)$  
$$
     \Lambda_\eff(\tau):=\left\{\lambda\in\Lambda(\tau)\,|\,
     \inner{\tau'}{\lambda}\geq 0\text{ for all }\tau'\in C(\tau)
     \right\},
$$
where $C(\tau)$ denotes the chamber of $\tau$. 
Note that 
$$
        \la\tau',\lambda\ra>0\quad\text{ for }
        \lambda\in\Lambda_\eff(\tau)\setminus\{0\},\ \tau'\in
        C(\tau).
$$
Denote by $\Dd_\eff(\tau)\subset\Z^n$ the cone
$$
     \Dd_\eff(\tau) := \left\{
     (\inner{\w_1}{\lambda},\dots,\inner{\w_n}{\lambda})\,|\,
     \lambda\in\Lambda_\eff(\tau)
     \right\}.
$$
Note that the map $\Lambda_\eff(\tau)\to\Dd_\eff(\tau)$
is a bijection.  We denote the inverse by 
$\Dd_\eff(\tau)\to\Lambda_\eff(\tau):d\mapsto\lambda_d$. 
We emphasize that $\Dd_\eff(\tau)$ is not 
necessarily contained in the positive quadrant of $\Z^n$. 

Let us now consider the vector
\begin{equation}\label{eq:monotone}
      \tau := \sum_{\nu=1}^n\w_\nu.
\end{equation}
(We still assume that $T$ acts freely on $\mu^{-1}(\tau)$.)
Then $\bar M$ is a monotone symplectic manifold
(see Lemma~\ref{le:monotone}).  The genus zero 
Gromov-Witten invariants of $\bar M$ with fixed marked points 
in a homology class $\bar\lambda\in H_2(\bar M)$ 
are denoted by 
$$
     \GW^{\bar M}_{\bar\lambda}:
     H^*(\bar M)\times\cdots\times H^*(\bar M)\to\Z. 
$$
The number of arguments will in each case be clear from the 
context. For an $n$-tuple $\ell=(\ell_1,\dots,\ell_n)$
of nonnegative integers and a cohomology class 
$\bar\alpha\in H^*(\bar M)$ we abbreviate
$$
     \GW^{\bar M}_{\bar\lambda}
     (\bar\w^{*\ell},\bar\alpha)
     := \GW^{\bar M}_{\bar\lambda}
     (\bar\w_1,\dots,\bar\w_1,\dots,
      \bar\w_n,\dots,\bar\w_n,\bar\alpha),
$$
where each argument $\bar\w_\nu$ occurs $\ell_\nu$ times.
Since the Gromov--Witten invariants are invariant under 
symplectic deformation, we have 
$$
\lambda\in\Lambda(\tau)\setminus\Lambda_\eff(\tau)
\qquad\IMP\qquad
\GW^{\bar M}_{\bar\lambda}\equiv0.
$$

\begin{lemma}\label{le:qcoh}
Let $\tau:=\sum_{\nu=1}^n\w_\nu$,
suppose that $T$ acts freely on $\mu^{-1}(\tau)$,
and assume that the minimal Chern number $N$
of $\bar M$ is greater that one. 
Then for every $d\in\Dd_\eff(\tau)$, 
every $\lambda\in\Lambda(\tau)$, and every 
$\bar\alpha\in H^*(\bar M)$,
we have
\begin{equation}\label{eq:qrel}
     \GW^{\bar M}_{\bar\lambda}
     (\bar\w^{*d^+},\bar\alpha)
     = \GW^{\bar M}_{\bar\lambda-\bar\lambda_d}
     (\bar\w^{*d^-},\bar\alpha),
\end{equation}
where the $n$-tuples $d^+$ and $d^-$ are defined by 
$$
     d_\nu^+ := \left\{\begin{array}{rl}
     d_\nu,&\mbox{if }d_\nu>0,\\
     0,&\mbox{if }d_\nu\le 0,
     \end{array}\right.\qquad
     d_\nu^- := \left\{\begin{array}{rl}
     -d_\nu,&\mbox{if }d_\nu<0,\\
     0,&\mbox{if }d_\nu\ge 0.
     \end{array}\right.\qquad
$$
\end{lemma}

\begin{proof}
Let $\bar\alpha:=\bar\w^\ell$ for an $n$-tuple $\ell$ of nonnegative
integers satisfying 
$$
        |\ell|=n-k+\la\tau,\lambda\ra-|d^+|.
$$
By Theorem~\ref{thm:invariants}~(iii), with
$\lambda$ replaced by $\lambda-\lambda_d$,
$\lambda'=\lambda_d$, and $\ell$ replaced by 
$\ell+d^-$, we have
\begin{equation}\label{eq:dpm}
     \Phi^{\rho,\tau}_\lambda(\w^{d^++\ell})
     = \Phi^{\rho,\tau}_{\lambda-\lambda_d}(\w^{d^-+\ell}).
\end{equation}
Hence it follows from~\cite[Theorem~A]{GS}
and the fact that $N>1$ (see equation~(\ref{eq:GS})
in the introduction) that
\begin{equation}\label{eq:GW1}
      \GW^{\bar M}_{\bar\lambda}(\bar\w^{*(d^++\ell)}) 
      = \GW^{\bar M}_{\bar\lambda-\bar\lambda_d}(\bar\w^{*(d^-+\ell)}).
\end{equation}
Now the gluing theorem for the Gromov--Witten invariants 
with fixed marked points (see~\cite{MS}) asserts that 
\begin{eqnarray*}
     \GW^{\bar M}_{\bar\lambda}(\bar\w^{*d^+},\bar\w^\ell) 
&= &
     \GW^{\bar M}_{\bar\lambda}(\bar\w^{*(d^++\ell)}) \\
&&-
     \sum_i\sum_{\lambda'\ne0}
     \GW^{\bar M}_{\bar\lambda-\bar\lambda'}(\bar\w^{*d^+},\bar e_i)
     \GW^{\bar M}_{\bar\lambda'}(\bar e_i^*,\bar\w^{*\ell}),
\end{eqnarray*}
where the second sum is over all lattice vectors 
$\lambda'\in\Lambda_\eff(\tau)\setminus\{0\}$. 
Hence, by~(\ref{eq:GW1}), 
\begin{eqnarray}\label{eq:GW2}
&&
     \GW^{\bar M}_{\bar\lambda}(\bar\w^{*d^+},\bar\w^\ell)
     - \GW^{\bar M}_{\bar\lambda-\bar\lambda_d}(\bar\w^{*d^-},\bar\w^\ell) \\
&&=
     \sum_i\sum_{\lambda'\ne0}
     \GW^{\bar M}_{\bar\lambda'}(\bar e_i^*,\bar\w^{*\ell})
     \left(
     \GW^{\bar M}_{\bar\lambda-\bar\lambda'-\bar\lambda_d}
     (\bar\w^{*d^-},\bar e_i)
     -
     \GW^{\bar M}_{\bar\lambda-\bar\lambda'}
     (\bar\w^{*d^+},\bar e_i)
     \right).
\nonumber
\end{eqnarray}
Note that in each summand on the right we have 
\begin{eqnarray*}
     \frac12\deg(e_i)
&= &
    n-k+\inner{\tau}{\lambda-\lambda'}-|d^+| \\
&< &
    n-k+\inner{\tau}{\lambda}-|d^+| \\
&= &
    |\ell|. 
\end{eqnarray*}
Hence the assertion follows from~(\ref{eq:GW2}) by induction 
over $|\ell|$.
\end{proof}

\noindent{\bf Remark.}
Let $\tau,\bar M$ be as in Lemma~\ref{le:qcoh} and
$\lambda\in\Lambda_\eff(\tau)$ such that 
$\inner{\w_\nu}{\lambda}\ge 0$ for every $\nu$. 
Then it follows from Lemma~\ref{le:qcoh} with 
$d_\nu:=\inner{\w_\nu}{\lambda}$ and 
$\bar\alpha:=\PD(\point)$ 
that $\GW^{\bar M}_{\bar\lambda}\ne0$. 
Hence the homology class $\bar\lambda\in H_2(\bar M)$ 
can be represented by a holomorphic stable map of genus zero.

\medskip

As in the introduction, let $\Rr$ be any graded commutative algebra
(over the reals) with unit which is equipped with a homomorphism
$$
     \Lambda_\eff(\tau)\to\Rr:\lambda\mapsto q^\lambda
$$
from the additive semigroup $\Lambda_\eff(\tau)$ to the 
multiplicative semigroup $\Rr$ such that 
$$
     \deg(q^\lambda) = 2\inner{\tau}{\lambda}. 
$$
The most important example is the ring 
$$
     \Rr=\R[q_1,\dots,q_k,q_1^{-1},\dots,q_k^{-1}]
$$
of polynomials with real coefficients in the variables 
$q_j$ and $q_j^{-1}$. To obtain the homomorphism choose
a basis $e_1,\dots,e_k$ of $\Lambda$, define the grading by 
$
\deg(q_j)=\sum_{\nu=1}^n\inner{\w_\nu}{e_j},
$
and the map $\lambda\mapsto q^\lambda$ by 
$$
q^\lambda:=\prod_{j=1}^kq_j^{\lambda_j},\qquad
\lambda=\sum_{j=1}^k\lambda_je_j.
$$ 
With a more careful choice of the basis 
one can take $\Rr=\R[q_1,\dots,q_k]$.
Other possibilities are the polynomial ring $\Rr=\R[q]$ 
in one variable, the ring of polynomials in $q$ and $q^{-1}$,
or the ring of Laurent series in $q$.  In these cases
one can choose $q$ to have degree two and define 
$q^\lambda:=q^{\inner{\tau}{\lambda}}$.
The simplest example is $\Rr=\R$ with the constant map
$\lambda\mapsto q^\lambda:=1$, but then the grading has to
be reduced modulo $2N$, where $N$ is the minimal Chern number.  

Given a graded algebra $\Rr$ as above define the 
quantum cohomology ring $\QH^*(\bar M;\Rr)$ as the 
tensor product 
$$
     \QH^*(\bar M;\Rr) := H^*(\bar M;\R)\otimes\Rr
$$
(of vector spaces over the reals).  Thus an element of 
$\QH^m(\bar M;\Rr)$ is a finite sum   
$
     \bar\alpha = \sum_{r\in\Rr}
     \bar\alpha_r r 
$
such that $\deg(\bar\alpha_r)+\deg(r)=m$ for all $r$. 
The ring structure is defined by 
$$
     \bar\alpha'*\bar\alpha''
     := \sum_i\sum_{\lambda\in\Lambda_\eff(\tau)}
        \sum_{r',r''}
        \GW^{\bar M}_{\bar\lambda}
        (\bar\alpha'_{r'},\bar\alpha''_{r''},\bar e_i)
        \bar e_i^* r'r''q^\lambda,
$$
where the $\bar e_i$ form a basis of $H^*(\bar M)$ and the 
$\bar e_i^*$ denote the dual basis with respect to 
the cup product pairing (see~\cite{MS}). 

\begin{corollary}\label{cor:qcoh}
Let $\tau:=\sum_{\nu=1}^n\w_\nu$,
suppose that $T$ acts freely on $\mu^{-1}(\tau)$,
and assume that the minimal Chern number $N$
of $\bar M$ is greater that one. Then
$$
     \bar\w^{*d^+} = \bar\w^{*d^-}q^{\lambda_d}
$$
for every $d\in\Dd_\eff(\tau)$. 
\end{corollary}

\begin{proof}
By the gluing theorem for the Gromov--Witten 
invariants (see~\cite{MS}), we have
\begin{eqnarray*}
     \bar\w^{*d^+} 
&=& 
     \sum_{i,\lambda}\GW^{\bar M}_{\bar\lambda}
     (\bar\w^{*d^+},\bar e_i)\bar e_i^*q^\lambda \\
&=& 
     \sum_{i,\lambda}\GW^{\bar M}_{\bar\lambda-\bar\lambda_d}
     (\bar\w^{*d^-},\bar e_i)\bar e_i^*q^\lambda \\
&=& 
     \bar\w^{*d^-}q^{\lambda_d}.
\end{eqnarray*}
The second equality follows from Lemma~\ref{le:qcoh}.
\end{proof}

\begin{proof}[Proof of Theorem~\ref{thm:batyrev}.]
We prove that the homomorphism~(\ref{eq:uq}) is surjective.
Note that there is an obvious inclusion 
$H^*(\bar M)\to 
\QH^*(\bar M;\Rr):
\bar\alpha\mapsto\bar\alpha 1$,
where $1$ denotes the unit in $\Rr$. Throughout we identify
$H^*(\bar M)$ with its image in 
$\QH^*(\bar M;\Rr)$ 
under this 
homomorphism.  Since~(\ref{eq:uq}) is a homomorphism 
of $\Rr$-modules, it suffices to prove that every class 
in $H^*(\bar M)$ belongs to the image of~(\ref{eq:uq}).  
We prove this by induction over the degree.
If $\bar\alpha\in H^0(\bar M)$ then $\bar\alpha$ 
obviously belongs to the image of~(\ref{eq:uq}).
Hence let $\deg(\bar\alpha)=2\ell>0$ and assume, by induction, that 
every class in $H^*(\bar M)$ of degree 
less than $2\ell$
belongs to the image of~(\ref{eq:uq}).
By Kirwan's theorem, the class $\bar\alpha$
is a linear combination of classes of the form 
$\bar\w_{\nu_1}\cdots\bar\w_{\nu_\ell}$. 
Let $p(u_1,\dots,u_n)$ be the same linear combination 
of the 
polynomials
$u_{\nu_1}\cdots u_{\nu_\ell}$. 
Then the image of $p(u)$ under the homomorphism 
differs from $\bar\alpha$ by a class of the form 
$$
     \bar\beta=\sum_{\lambda\ne 0}\bar\beta_\lambda q^\lambda,\qquad
     \deg(\bar\beta_\lambda) = 2\ell-2\inner{\tau}{\lambda}<2\ell.
$$
Here the sum is over all $\lambda\in\Lambda_\eff(\tau)$
that satisfy $\inner{\tau}{\lambda}>0$.  
Hence, by the induction hypothesis, every $\bar\beta_\lambda$ 
in this sum belongs to the image of~(\ref{eq:uq}),
and so does the class $\bar\beta_\lambda q^\lambda$.  
Hence $\bar\beta$ belongs to the image of~(\ref{eq:uq}),
and so does $\bar\alpha$. 

Let $\Jj_0\subset\Rr[u_1,\dots,u_n]$ be the kernel 
of~(\ref{eq:uq}). Then the linear polynomial 
$\sum_\nu\eta_\nu u_\nu$ belongs to $\Jj_0$ whenever 
$\sum_\nu\eta_\nu\w_\nu=0$.  Moreover, 
by Corollary~\ref{cor:qcoh}, the polynomial 
$u^{d^+}-q^\lambda u^{d^-}$ belongs to $\Jj_0$ whenever
$\lambda\in\Lambda_\eff(\tau)$ and 
$d_\nu^\pm=\max\{\pm\inner{\w_\nu}{\lambda},0\}$.
Hence $\Jj\subset\Jj_0$.

We prove that $\Jj_0\subset\Jj$.  
Define the classes $\bar\alpha_{\ell,\lambda}\in H^*(\bar M)$,
for 
$n$-tuples $\ell$ of nonnegative integers and 
lattice vectors $\lambda\in\Lambda_\eff(\tau)$ with 
$0<\inner{\tau}{\lambda}\le|\ell|$, by 
$$
     \bar\w^{*\ell} =: \bar\w^\ell 
     + \sum_{\lambda\in\Lambda_\eff(\tau)\setminus\{0\}}
       \bar\alpha_{\ell,\lambda} q^\lambda,\qquad
     \deg(\bar\alpha_{\ell,\lambda}) = 2|\ell|-2\inner{\tau}{\lambda}.
$$
For $N\in\Z$ denote by $\Jj_0(N)$ the set of polynomials
$p\in\Jj_0$ of the form 
\begin{equation}\label{eq:p}
      p(u_1,\dots,u_n)=\sum_{|\ell|\le N} r_\ell u^\ell,
\end{equation}
where the sum is over all $n$-tuples 
$\ell=(\ell_1,\dots,\ell_n)$ of nonnegative integers
satisfying $|\ell|\le N$. We prove by induction on $N$ 
that $\Jj_0(N)\subset\Jj$.
For $N<0$ this is obvious because $\Jj_0(N)=\{0\}$.
Let $N\ge 0$ and assume by induction that $\Jj_0(N-1)\subset\Jj$.  
Let $p\in\Jj_0(N)$ be a polynomial of the form~(\ref{eq:p}). 
Since $p\in\Jj_0$ we have 
$$
     0 = \sum_\ell r_\ell\bar\w^{*\ell}
       = \sum_\ell \bar\w^\ell r_\ell 
       + \sum_\ell\sum_\lambda\bar\alpha_{\ell,\lambda}r_\ell q^\lambda.
$$
This identity splits up into 
$$
     0 = \sum_{|\ell|=j} r_\ell\bar\w^\ell 
         + \sum_{|\ell|>j}
           \sum_{\lambda\atop\inner{\tau}{\lambda}=|\ell|-j}
          \bar\alpha_{\ell,\lambda}r_\ell q^\lambda,\qquad
     j=0,\dots,N.
$$
Since $r_\ell=0$ for $|\ell|>N$, we have 
$$
     \sum_{|\ell|=N}r_\ell\bar\w^\ell =0.
$$
Choose a basis $\rho_1,\dots,\rho_m$ of the vector space 
$\SPAN\{r_\ell\,|\,|\ell|=N\}\subset\Rr$ and express each 
$r_\ell$ in this basis, i.e.
$$
     r_\ell = \sum_{i=1}^ma_{\ell i}\rho_i,\qquad
     a_{\ell i}\in\R,\qquad|\ell|=N.
$$
Then 
$$
     \sum_{|\ell|=N}a_{\ell i}\bar\w^\ell=0,\qquad
     i=1,\dots,m.
$$
This means that the polynomials 
$$
     p_{i0}(u_1,\dots,u_n):=\sum_{|\ell|=N}a_{\ell i}u^\ell,\qquad
     i=1,\dots,m,
$$
belong to the kernel $\Ii\subset\R[u_1,\dots,u_n]$ of the
homomorphism~(\ref{eq:u}) in Theorem~\ref{thm:coh}. 
Hence they can be expressed in the form
$$
     p_{i0} = \sum_jp_{ij0}f_j,
$$
where $f_j\in\R[u_1,\dots,u_n]$, and the $p_{ij0}$ 
are taken from the set of generators of $\Ii$ in 
Theorem~\ref{thm:ideal}. Thus each $p_{ij0}$ satisfies 
one of the following conditions.
\begin{description}
\item[(a)]
$p_{ij0}(u)=\sum_\nu\eta_\nu u_\nu$, 
where $\sum_\nu\eta_\nu\w_\nu=0$.
\item[(b)]
$p_{ij0}(u)=u_\nu$, where $\Delta_{\{\nu\}}=\emptyset$. 
\item[(c)]
$p_{ij0}(u)=u^{d^+}$, where $d\in\Dd_\eff(\tau)\setminus\{0\}$.
\end{description}
In cases~(a) and~(b) define $p_{ij}:=p_{ij0}\in\Jj$.
In the case~(c) it follows from the definition of $\Jj$ 
that there is a generator $p_{ij}\in\Jj$ of the form 
$$
      p_{ij}(u) = u^{d^+}-q^\lambda u^{d^-},\qquad
      p_{ij0}(u)=u^{d^+},\qquad 
      d\in\Dd_\eff(\tau)\setminus\{0\}.
$$
Define $\hat p\in\Jj$ by 
$$
      \hat p(u) := \sum_{i=1}^m\sum_j
      \rho_ip_{ij}(u)f_j(u).
$$
Since $\Jj\subset\Jj_0$ we have $p-\hat p\in\Jj_0$. 
Since 
$$
      \sum_{i=1}^m\sum_j\rho_ip_{ij0}(u)f_j(u)
      = \sum_{i=1}^m\rho_ip_{i0}(u)
      = \sum_{i=1}^m\sum_{|\ell|=N}\rho_ia_{\ell i}u^\ell
      = \sum_{|\ell|=N}r_\ell u^\ell,
$$
the leading terms cancel in $p-\hat p$ and hence 
$p-\hat p\in\Jj_0(N-1)\subset\Jj$.  Hence $p\in\Jj$. 
This completes the induction and the proof of the theorem. 
\end{proof}

\begin{example}\rm
This example shows that in the definition of the 
ideal $\Jj$ it may not suffice to consider vectors 
$\lambda\in\Lambda_\eff(\tau)$ such that the integers
$d_\nu:=\inner{\w_\nu}{\lambda}$ are all nonnegative.
Suppose the 2-torus $T=\T^2$ acts on $\C^5$ 
with weight vectors 
$$
    \w_1=(1,0),\qquad \w_2=(1,1),\qquad \w_3=\w_4=\w_5=(0,1).
$$
The symplectic quotient $\bar M$ at the parameter 
$
     \tau:=\w_1+\cdots+\w_5=(2,4)
$
is a smooth monotone toric $3$-fold with minimal 
Chern number $N=2$. The effective cone is given by 
$$
     \Lambda_\eff(\tau)=\left\{(\lambda_1,\lambda_2)\in\Z^2\,|\,
     \lambda_2\ge0,\,\lambda_1+\lambda_2\ge 0\right\}.
$$
It is the convex cone spanned by the vectors
$e:=(1,0)$ with $d=(1,1,0,0,0)$ and
$e'=(-1,1)$ with $d'=(-1,0,1,1,1)$. 
For the quantum cohomology let us choose the polynomial 
ring $\Rr:=\R[q_1,q_2]$, graded by $\deg(q_1)=\deg(q_2)=4$,
and the homomorphism 
$
     q^\lambda:=q_1^{\lambda_1+\lambda_2}q_2^{\lambda_2}.
$
Thus $q_1,q_2$ correspond to the generators $e,e'$ of
$\Lambda_\eff(\tau)$. Then the ideal
$\Jj\subset\R[u_1,\dots,u_5,q_1,q_2]$ 
is generated by the relations
\begin{equation}\label{eq:rel}
     u_3=u_4=u_5=u_2-u_1,\qquad
     u_1u_2=q_1,\qquad u_3u_4u_5=u_1q_2.
\end{equation}
If one considers only vectors $\lambda\in\Lambda_\eff(\tau)$ 
with nonnegative degrees $d_\nu:=\inner{\w_\nu}{\lambda}$
then one has to replace the last relation 
in~(\ref{eq:rel}) by $u_2u_3u_4u_5=q_1q_2$
and obtains a strictly smaller ideal. 
\end{example}


\appendix

\section{Equivariant trivialization}

\begin{proposition}\label{prop:triv}
Let $\G$ be a compact Lie group and $E\to X$
be a $\G$-equi\-var\-i\-ant complex vector bundle
over  compact smooth manifold $X$.
Then there exists a $\G$-equivariant complex vector bundle
$F\to X$ and a complex $\G$-representation $W$ such that
$E\oplus F$ is equivariantly isomorphic to $X\times W$.
\end{proposition}

\begin{lemma}\label{le:VW}
Let $\G$ be a compact Lie group, $\HG\subset\G$ be a
normal subgroup, and $V$ be a complex $\HG$-representation.
Then there exists a complex $\G$-representation $W$
and an injective $\HG$-equivariant homorphism
$\phi:V\to W$.
\end{lemma}

\begin{proof}
Consider the infinite dimensional vector space
$$
    \Ww := \left\{f:\G\to V\,|\,f(hg)=hf(g)\;
    \forall h\in\HG,\;\forall g\in\G\right\}.
$$
This space carries an action of $\G$ by
$$
    (g'f)(g) := f(gg')
$$
and the evaluation map $\Ww\to V:f\mapsto f(1)$
is $\HG$-equivariant and surjective.  To prove surjectivity
let $v\in V$ be given and let $f:\G\to V$ be any smooth
extension of the map $\HG\to V:h\mapsto hv$.  By averaging
the maps $g\mapsto h^{-1}f(hg)$ over $h\in\HG$
we can ensure that the extension is $\HG$-equivariant.
By Peter--Weyl's theorem, there exists a finite dimensional
$\G$-invariant subspace $W\subset\Ww$ such that the restriction
of the homomorphism $f\mapsto f(1)$ to $W$ is still surjective.
By~\cite[Remark~A.4.2]{MS}, the surjection $W\to V$
has an $\HG$-equivariant right inverse.
\end{proof}

\begin{proof}[Proof of Proposition~\ref{prop:triv}]
Let $x\in X$, $\HG\subset\G$ be the isotropy
subgroup of $x$, and $V:=E_x$. By the local slice theorem,
the restriction of $E$ to a suitable neighbourhood of
the $\G$-orbit of $x$ is equivariantly isomorphic to
the bundle
$$
     \frac{\G\times U\times V}{\HG}\to \G\times_\HG U,
$$
where $U$ is a neighbourhood of zero in the horizontal
tangent space at $x$ (i.e. in the orthogonal complement
of $T_x\G x$ with respect to some $\G$-invariant metric).
Let $\phi:V\to W$ be as in Lemma~\ref{le:VW}.
Then the map
$$
     \G\times U\times V\to \G\times U\times W:
     (g,u,v)\mapsto(g,u,g\phi(v))
$$
descends to a $\G$-equivariant injective bundle
homomorphism from $(\G\times U\times V)/\HG$
to $(\G\times_\HG U)\times W$ (where $\G$ acts diagonally).
This construction gives rise to a $\G$-invariant open cover
$\{U_\alpha\}_\alpha$ of $X$ and a collection of
$\G$-equivariant injective bundle homomorphisms
$\phi_\alpha:E|_{U_\alpha}\to U_\alpha\times W_\alpha$.
Let $\rho_\alpha:X\to[0,1]$ be a $\G$-invariant partition
of unity subordinate to the cover $\{U_\alpha\}_\alpha$
and denote $W:=\bigoplus_\alpha W_\alpha$.  Then the
homomorphism
$
    E\to X\times W:
    (x,e)\mapsto (x,\{\rho_\alpha(x)\phi_\alpha(x)e\}_\alpha)
$
is the required $\G$-equivariant embedding.
\end{proof}


\section{Convex polytopes}\label{app:convex}

In this section we recall some well-known facts about convex polytopes
(see e.g.~\cite{ODA}). Let $\Delta$ be a compact convex polytope in
the dual space $V^*$ of a finite dimensional vector space $V$. 
We denote elements of $V$ by $v,w$ and elements of $V^*$ 
by $\xi,\eta$. Define the {\em support function} 
$\phi:V\to\R$ of $\Delta$ by
$$
     \phi(v) := \inf_{\xi\in\Delta}\la\xi,v\ra.
$$ 
The following properties of $\phi$ are obvious from the definition.
\begin{description}
\item[(P1)] 
$\phi(tv)=t\phi(v)$ for $t\geq 0$.
\item[(P1)] 
$\phi$ is concave, i.e.~$\phi(v+w)\geq\phi(v)+\phi(w)$.
\item[(P3)] 
$\Delta$ can be recovered from $\phi$ as the intersection of
half spaces
$$
     \Delta = \bigcap_{v\in V}\{\xi\in V^*\;|\;\la\xi,v\ra\geq\phi(v)\}. 
$$
\end{description}
Let $F$ be a face of $\Delta$. Pick an interior point $p$ of $F$
and define the {\em dual cone to $F$}
by 
$$
     \check{F} := \{v\in V\,|\,\la\xi-p,v\ra\geq 0
     \text{ for all }\xi\in\Delta\}. 
$$ 
If $q$ is another interior point of $F$ and $\xi\in\Delta$ then
$q+t(\xi-p)\in\Delta$ for $t>0$ sufficiently small. Hence
$\la\xi-p,v\ra\geq 0$ iff $\la q+t(\xi-p)-q,v\ra\geq 0$. This shows
that the definition of $\check{F}$ does not depend on the point
$p$. Moreover, the condition 
$
     \la\xi-p,v\ra\ge0\text{ for all }\xi\in\Delta
$ 
can be rewritten as $\phi(v)\geq\la p,v\ra$, or equivalently
$\phi(v)=\la p,v\ra$ since $p\in\Delta$. So $\check{F}$
can be written in the equivalent forms
\begin{eqnarray*}
\check{F} 
&= &
\{v\in V\,|\,\la\xi-p,v\ra\geq 0\text{ for all }\xi\in\Delta,p\in F\} \\
&= &
\{v\in V\,|\,\la p,v\ra=\phi(v)\text{ for all }p\in F\}.
\end{eqnarray*}
The following properties are obvious from these descriptions 
of $\check{F}$. 
\begin{description}
\item[(F1)] 
$\check{F}$ is a convex polyhedral cone.
\item[(F2)] 
The restriction of $\phi$ to $\check{F}$ is the linear function
$\phi(v)=\la p,v\ra$ for any $p\in F$.
\item[(F3)]
$\check{F}$ is perpendicular to $F$ and 
$\dim\check{F}=\codim F$.
\item[(F4)]
If $H_1,\dots,H_\ell\subset V^*$ are the supporting hyperplanes 
for $\Delta$ meeting at $F$, then $\check{F}$ is the cone generated 
by inward pointing normal vectors $v_1,\dots,v_\ell$ to the hyperplanes.
\item[(F5)]
If $G$ is a subface of $F$ then $\check{F}$ is a subcone of
$\check{G}$. 
\item[(F6)]
The union of the cones $\check{p}$ dual to vertices $p$ of $\Delta$
is the whole space $V$.
\end{description}
The collection $\Sigma$ of the cones $\check{F}$ dual to nonempty
faces of $\Delta$ is called the {\em fan dual to $\Delta$}
(see~\cite{AUDIN} for the general definition of a fan). 


\section{The cohomology of symplectic quotients}\label{app:coh}

Let $T$ be a $k$-dimensional torus and
$\rho=(\rho_1,\dots,\rho_n):T\to\T^n:=(S^1)^n$ be a diagonal  
homomorphism with
$$
     \rho_\nu(\exp(\xi)) = e^{-2\pi i\la\w_\nu,\xi\ra}
$$
for $\xi\in\tt:=\Lie(T)$. Here the $\w_\nu$ are elements 
of the dual lattice $\Lambda^*\subset\tt^*$ as in the 
introduction.  We identify the Lie algebra of 
$\T^n$ with $\R^n$ via the map $\eta\mapsto i\eta/2\pi$ 
so that the integer lattice corresponds to $\Z^n\subset\R^n$. 
In this identification the linearization of $\rho$ is the map
$
        \dot\rho:\tt\to\R^n
$
given by 
$$
     \dot\rho(\xi) = \bigl(\la\w_1,\xi\ra,\dots,\la\w_n,\xi\ra\bigr). 
$$
Consider the quotient torus
$$
    \bar T := \T^n/\rho(T).
$$
Its Lie algebra is the quotient space 
$$
    \bar\tt:=\R^n/\dot\rho(\tt)
$$
and the dual space of $\bar\tt$ can be identified with the 
subspace 
\begin{equation}\label{eq:tbar}
    \bar\tt^* := \left\{\eta\in(\R^n)^*\,\Big|\,
    \sum_{\nu=1}^n\eta_\nu\w_\nu=0\right\}.
\end{equation}
The canonical action of $\T^n$ on $\C^n$ induces an action of $T$ with
moment map $\mu:\C^n\to\tt^*$ 
given by~(\ref{eq:mu}). We assume throughout that $\mu$ is proper
and that the action is effective (i.e. the weight vectors $\w_\nu$
span $\tt^*$). Let $\tau\in\tt^*$ be a regular value 
of $\mu$. Then the torus $\bar T$ acts on the 
symplectic quotient
$$
    \bar M := \C^n\dslash T(\tau) = \mu^{-1}(\tau)/T.
$$
A moment map $\bar\mu:\bar M\to\bar\tt^*$ for this action 
is given by the formula
$$
    \bar\mu(x) := \left(\begin{array}{c}
     \pi|x_1|^2+\zeta_1 \\ \vdots \\ \pi|x_n|^2+\zeta_n
    \end{array}\right),
$$
where $\zeta=(\zeta_1,\dots,\zeta_n)\in(\R^n)^*$
is chosen such that 
$$
    \sum_{\nu=1}^n\zeta_\nu\w_\nu = -\tau.
$$
The image of $\bar\mu$ is the convex polyhedron
\begin{equation}\label{eq:Delta}
    \Delta := \bar\mu(\bar M)
    = \left\{\eta\in(\R^n)^*\,\Big|\,\sum_{\nu=1}^n\eta_\nu\w_\nu=0,\,
      \eta_\nu\geq\zeta_\nu\right\}.
\end{equation}
Each subset $I\subset\{1,\dots,n\}$ determines a (possibly empty) face
$$
        \Delta_I := \left\{\eta\in\Delta\,|\,
     \eta_\nu=\zeta_\nu\text{ for }\nu\in I\right\}.
$$
Recall that $C(I)$ denotes the cone spanned by the vectors $\w_\nu$,
$\nu\in I$. The next lemma shows that if $\tau$ is a regular value of $\mu$,
then the intersection of any $j$ codimension-$1$ faces of $\Delta$
is either empty or has codimension $j$. 

\begin{lemma}\label{le:regular}
Assume that $\tau$ is a regular value of $\mu$ and
let $I\subset\{1,\dots,n\}$.  

\smallskip
\noindent{\bf (i)}
The set $\Delta_I$ is either empty or has codimension $|I|$. 

\smallskip
\noindent{\bf (ii)}
$
        \Delta_I=\emptyset \Longleftrightarrow 
        \tau\notin C(\{1,\dots,n\}\setminus I). 
$
\end{lemma}

\begin{proof}
We prove~(i). Assume $\Delta_I\ne\emptyset$ and let 
$J:=\{1,\dots,n\}\setminus I$. Then, by the definition 
of $\bar\mu$, there is a $y\in\C^J$ such that  
$$
     \mu(y)=\pi\sum_{\nu\in J}|y_\nu|^2\w_\nu=\tau.
$$
Since $\tau$ is a regular value of $\mu$, there exist
indices $\nu_1,\dots,\nu_k\in J$ such that the vectors
$\w_{\nu_1},\dots,\w_{\nu_k}$ are linearly independent
and $y_{\nu_j}\ne 0$ for every $j$. 
We claim that there is a vector 
$x\in\C^J$
such that
$$
     \mu(x) = \tau,\qquad x_\nu\ne 0\quad\mbox{for all}\;\;\nu\in J.
$$
To see this choose 
$x_\nu$ for $\nu\in J\setminus\{\nu_1,\dots,\nu_k\}$
such that $|x_\nu|^2=|y_\nu|^2+\eps$
and choose $x_{\nu_j}$ such that 
$$
     \sum_{j=1}^k(|x_{\nu_j}|^2-|y_{\nu_j}|^2)\w_{\nu_j}
     + \eps\sum_{\nu\in J\setminus\{\nu_1,\dots,\nu_k\}}\w_\nu
     = 0.
$$
Then $\mu(x)=\mu(y)=\tau$ and, for $\eps>0$ sufficiently small, 
we get $x_\nu\ne 0$ for all $\nu\in J$. 

The differential $d\bar\mu(x):T_x(\bar M\cap\C^J)\to\bar\tt^*$ is
given by
$$
        d\bar\mu(x)v = (2\pi\la x_\nu,v_\nu\ra)_{\nu\in J},
$$
where $v\in\C^J$ satisfies 
$$
        d\mu(x)v = 2\pi\sum_{\nu\in J}\la x_\nu,v_\nu\ra \w_\nu = 0.
$$
Since $x_\nu\neq 0$ for $\nu\in J$, this shows that the image of
$T_x(\bar M\cap\C^J)$ under $d\bar\mu(x)$ equals
$\{\eta\in(\R^J)^*\;|\;\sum_{\nu\in J}\eta_\nu\w_\nu=0\}$. This space,
and therefore $\Delta_I$, has dimension $|J|-k = n-k-|I|$.

We prove~(ii). If $\Delta_I\neq\emptyset$ there exists an
$\eta\in(\R^n)^*$ such that 
$
     \sum_{\nu=1}^n\eta_\nu\w_\nu=0,
$
$
     \eta_\nu\geq\zeta_\nu
$
for all $\nu$, and 
$
    \eta_\nu=\zeta_\nu
$
for $\nu\in I$.  Hence
$$
     \tau = \sum_{\nu=1}^n-\zeta_\nu w_\nu = \sum_{\nu\notin
     I}(\eta_\nu-\zeta_\nu)\w_\nu 
     \in C(\{1,\dots,n\}\setminus I). 
$$
The converse follows by reversing the argument.
\end{proof}

\medskip
\noindent{\bf Standing assumption.}
{\it In the remainder of this appendix we assume 
that $T$ acts freely on $\mu^{-1}(\tau)$.}

\medskip
\noindent
Denote by $\bar\w_\nu\in H^2(\bar M;\R)$ the image of $\w_\nu$
under the homomorphism 
$\Lambda^*\cong H^2(\BT;\Z)\to H^*(\bar M;\R)$.

\begin{lemma}\label{le:cup}
For every $J\subset\{1,\dots,n\}$ the following holds.

\smallskip
\noindent{\bf (i)}
If $\tau\notin C(J)$ then 
$
      \prod_{\nu\notin J}\bar\w_\nu = 0.
$

\smallskip
\noindent{\bf (i)}
If $\tau\in C(J)$ and $|J|=k$ then 
$
      \prod_{\nu\notin J}\bar\w_\nu
      = \PD(\point).
$
\end{lemma}

\begin{proof}
$
     \bar\w_\nu
$
is the first Chern class of the line bundle 
$
     \bar L_\nu := \mu^{-1}(\tau)\times_{\rho_\nu}\C.
$
Hence the zero set of the holomorphic section 
$\bar M\to \bar L_\nu:[x]\mapsto[x,x_\nu]$
is Poincar\'e dual to $\bar\w_\nu$.  Denote this zero
set by 
$
     \bar W_\nu := \left\{[x]\in\bar M\,|\,x_\nu=0\right\}.
$
This is a (possibly empty) complex submanifold of 
$\bar M$ of complex codimension one.  Moreover,
$$
     \tau\notin C(J)\qquad\IMP\qquad
     \bigcap_{\nu\notin J}\bar W_\nu=\emptyset
$$
This proves~(i).  If $\tau\in C(J)$ and $|J|=k$, then the 
submanifolds $\bar\W_\nu$ for $\nu\notin J$ intersect 
transversally in a single point. This proves~(ii). 
\end{proof}

\begin{lemma}\label{le:monotone}
{\bf (i) }The Chern classes of $T\bar M$ are given by 
$$
     c_j(T\bar M) = \sum_{\nu_1<\cdots<\nu_j}
     \bar\w_{\nu_1}\cdots\bar\w_{\nu_j},\qquad
     j=1,\dots,n-k.
$$
{\bf (ii)}
The cohomology class of the symplectic form 
$\bar\om\in\Om^2(\bar M)$ is 
$$
     [\bar\om] = \bar\tau. 
$$
\end{lemma}

\begin{proof}
We prove~(i). Consider the Whitney sum
$
     \mu^{-1}(\tau)\times\C^n=E\oplus F,
$
where the complex vector bundles $E\to\mu^{-1}(\tau)$
and $F\to\mu^{-1}(\tau)$ are defined by 
$$
     E_x := \left\{v\in\C^n\,|\,d\mu(x)v=d\mu(x)iv=0\right\},
$$
$$
     F_x := \left\{v\in\C^n\,|\,\exists\xi,\eta\in\tt\;\forall\nu:
     v_\nu = (\inner{\w_\nu}{\xi}+i\inner{\w_\nu}{\eta})x_\nu\right\}.
$$
Then the bundle $F$ admits a $T$-equivariant complex 
trivialization and the quotient bundle $E/T\to\mu^{-1}(\tau)/T$
is isomorphic to the tangent bundle of $\bar M$.  Hence 
$$
     c_j(T\bar M)
     = c_j(\mu^{-1}(\tau)\times_T\C^n)
     = \sum_{\nu_1<\cdots<\nu_j}
     \bar\w_{\nu_1}\cdots\bar\w_{\nu_j}.
$$

We prove~(ii). Denote
$$
     \lambda_0 
     := \frac{1}{2i}\sum_{\nu=1}^n(\bar x_\nu dx_\nu
        - x_\nu d\bar x_\nu)
     \in\Om^1(\C^n),\qquad
     \om_0 := d\lambda_0,
$$
and let $\mu_0:\C^n\to\R^n$ be the moment map given by 
$$
     \mu_0(x) := \pi\left(|x_1|^2,\dots,|x_n|^2\right).
$$
Since $d_{\T^n}\lambda_0=\om_0-\mu_0$, the equivariant
cohomology class $[\om_0-\mu_0]\in H^2_{\T^n}(\C^n)$ is
trivial. Pulling back under the homomorphism
$H^2_{\T^n}(\C^n)\to H^2_T(\C^n)$ induced by 
$\rho$ yields $0=[\om_0-\mu]\in H^2_T(\C^n)$. 
Restriction to $\mu^{-1}(\tau)$ yields
$0=[\iota^*\om_0-\iota^*\tau]\in H^2_T\bigl(\mu^{-1}(\tau)\bigr)$,
where $\iota:\mu^{-1}(\tau)\to\C^n$ is the inclusion. 
Now the result follows by passing to the quotient.
\end{proof}

\begin{theorem}[\cite{GuSt}]\label{thm:coh}
The ring homomorphism 
\begin{equation}\label{eq:u}
\R[u_1,\dots,u_n]\to H^*(\bar M;\R):
p(u_1,\dots,u_n)\mapsto p(\bar\w_1,\dots,\bar\w_n)
\end{equation}
induces an isomorphism
$$
      H^*(\bar M;\R)\cong\R[u_1,\dots,u_n]/\Ii,
$$
where the ideal $\Ii\subset\R[u_1,\dots,u_n]$
is generated by the relations
\begin{equation}\label{eq:rel0}
     \sum_{\nu=1}^n\eta_\nu\w_\nu=0\qquad\IMP\qquad
     \sum_{\nu=1}^n\eta_\nu u_\nu=0,
\end{equation}
\begin{equation}\label{eq:rel1}
     I\subset\{1,\dots,n\},\quad
     \Delta_I = \emptyset\qquad\IMP\qquad
     \prod_{\nu\in I}u_\nu=0.
\end{equation}
\end{theorem}

\noindent{\bf Remark.}
Theorem~\ref{thm:coh} continues to hold with coefficients in $\Z$. 
That the homomorphism~(\ref{eq:u}) is surjective follows from
Kirwan's theorem, and that the ideal $\Ii$ is contained 
in the kernel of~(\ref{eq:u}) is an easy consequence of
Lemma~\ref{le:cup}. The nontrivial part of the proof is to show
that the kernel is contained in $\Ii$. 

\begin{theorem}[\cite{BATYREV}]\label{thm:ideal}
Assume $\mu^{-1}(\tau)\ne\emptyset$ 
and $T$ acts freely on $\mu^{-1}(\tau)$.  
Then the ideal $\Ii$ is generated by the linear 
relations~(\ref{eq:rel0}), the linear monomials $u_\nu$
for $\Delta_{\{\nu\}}=\emptyset$, and the monomials 
$u^{d^+}$ for $d\in\Dd_\eff(\tau)\setminus\{0\}$.  
\end{theorem}

For the sake of completeness, we present a
somewhat more elaborated version of the proof 
given in~\cite{BATYREV}. We need some preparation. 
Denote by $\bar\Lambda\subset\bar\tt$ 
and $\bar\Lambda^*\subset\bar\tt^*$ the integer lattices. 
Thus $\bar\Lambda$ is the image of $\Z^n$ under the projection
$\R^n\to\bar\tt$ and $\bar\Lambda^*=\bar\tt^*\cap(\Z^n)^*$. 
For $\nu=1,\dots,n$ let $\bar e_\nu\in\bar\Lambda$ be the image of 
the basis vector $e_\nu=(0,\dots,0,1,0,\dots,0)\in\Z^n$
under the projection $\R^n\to\bar\tt$.

\begin{lemma}\label{le:lattice}
Suppose that $T$ acts freely on $\mu^{-1}(\tau)$.

\smallskip
\noindent{\bf (i) }Let $J\subset\{1,\dots,n\}$ satisfy $|J|=n-k$ and
$\Delta_J\ne\emptyset$. Then the vectors $\{\bar
e_j\,|\,j\in J\}$ form an integer basis of $\bar\Lambda$. 

\smallskip
\noindent{\bf (ii) }Let $d_\nu\in\Z$ satisfy $\sum_{\nu=1}^nd_\nu\bar e_\nu=0$.
Then there exists a vector $\lambda\in\Lambda$ 
such that $d_\nu=\inner{\w_\nu}{\lambda}$ for every $\nu$.
\end{lemma}

\begin{proof}
We prove (i). Assume $|J|=n-k$ and $\Delta_J\ne\emptyset$.
Since $T$ acts freely on $\mu^{-1}(\tau)$, and 
$\tau\in C(\{1,\dots,n\}\setminus J)$, the vectors
$\{\w_\nu\,|\,\nu\notin J\}$ form an integer basis of 
$\Lambda^*$. 
Hence, for every $v\in\Z^n$, there exists a unique vector 
$\lambda\in\Lambda$ such that $v_\nu=\inner{\w_\nu}{\lambda}$ for
$\nu\notin J$. 
This implies that the image $\bar v\in\bar\tt$ of $v$ 
under the projection $\R^n\to\bar\tt$ satisfies
$$
      \bar v = \sum_{j\in J}
      (v_j-\inner{\w_j}{\lambda})\bar e_j.
$$
Hence the vectors $\{\bar e_j\,|\,j\in J\}$ span the 
integer lattice $\bar\Lambda$ as claimed.

We prove (ii). By definition of the projection $\R^n\to\bar\tt$, there
exists a vector $\xi\in\tt$ such that $d_\nu=\inner{\w_\nu}{\xi}$
for every $\nu$.  Now let $J\subset\{1,\dots,n\}$ be any index set
such that $|J|=n-k$ and $\Delta_J\ne\emptyset$.  Then 
the argument in the proof of (i) shows that there exists 
a lattice vector $\lambda\in\Lambda$ such that 
$d_\nu=\inner{\w_\nu}{\lambda}$ for $\nu\notin J$. 
Hence $\inner{\w_\nu}{\xi-\lambda}=0$ for $\nu\notin J$.
Since the vectors $\{\w_\nu\,|\,\nu\notin J\}$ form a basis
of $\tt^*$ we deduce that $\xi=\lambda$ and hence
$d_\nu=\inner{\w_\nu}{\lambda}$ for every $\nu$. 
\end{proof}

\begin{proof}[Proof of Theorem~\ref{thm:ideal}]
Let $\Ii_0\subset\R[u_1,\dots,u_n]$ be the ideal 
generated by the linear polynomials 
$\sum_{\nu=1}^n\eta_\nu u_\nu$, where 
$\sum_{\nu=1}^n\eta_\nu\w_\nu=0$,
the monomials $u_\nu$, where $\Delta_{\{\nu\}}=\emptyset$,
and the monomials $u^{d^+}$ for $d\in\Dd_\eff(\tau)\setminus\{0\}$.  

We prove that $\Ii_0\subset\Ii$. 
We must show that $u^{d^+}\in\Ii$ for every
$d\in\Dd_\eff(\tau)\setminus\{0\}$. 
We prove a stronger statement: 
{\it If $\lambda\in\Lambda$ satisfies $\inner{\tau}{\lambda}>0$ 
and $d_\nu:=\inner{\w_\nu}{\lambda}$ then $u^{d^+}\in\Ii$.}
To see this, consider the set 
$
    I:=\{\nu\;|\;d_\nu>0\}.
$
We claim that $\Delta_I=\emptyset$. 
Otherwise, by Lemma~\ref{le:regular}, there
would exist numbers $\eta_\nu\geq 0$ such that 
$\tau=\sum_{\nu\notin I}\eta_\nu\w_\nu$. But then 
$$
    0 < \la\tau,\lambda\ra 
    = \sum_{\nu\notin I}\eta_\nu d_\nu 
    \le 0,
$$
a contradiction. Since $\Delta_I=\emptyset$, the monomial 
$\prod_{\nu\in I}u_\nu$ belongs to the ideal $\Ii$. 
But $u^{d^+}$ is a multiple of $\prod_{\nu\in I}u_\nu$ 
and hence also belongs to $\Ii$. 

We prove that $\Ii\subset\Ii_0$.  Consider the moment
polytope $\Delta\subset\bar\tt^*$ defined by~(\ref{eq:Delta}). 
The faces of $\Delta$ are subsets of the 
form $\Delta_I$ for $I\subset\{1,\dots,n\}$
such that $\tau\in C(\{1,\dots,n\}\setminus I)$.  
The vectors $\{\bar e_i\,|\,i\in I\}$ are the inward 
pointing normal vectors to the supporting hyperplanes of 
$\Delta$ meeting at the face $\Delta_I$. Hence, by property~(F4) 
of the dual cones (see Appendix~\ref{app:convex}), the dual cone 
of $\Delta_I$ is given by 
$$
     \check{\Delta}_I
     = \left\{\sum_{i\in I}c_i\bar e_i\,\Big|\,c_i\ge 0\right\}. 
$$
By Lemma~\ref{le:regular}, the codimension of the face $\Delta_I$
equals $|I|$. In particular, the vertices of $\Delta$ are subsets
$\Delta_J$ where $|J|=n-k$ and $\Delta_J\ne\emptyset$. 

Now let $I\subset\{1,\dots,n\}$ such that $\Delta_I=\emptyset$. 
We must prove that the monomial $\prod_{\nu\in I}u_\nu$ 
belongs to $\Ii_0$. Shrinking the set $I$, if necessary, we may
assume without loss of generality that $\Delta_{I'}\neq\emptyset$ 
for every proper subset $I'\subsetneq I$.  Since $\mu(\tau)\ne\emptyset$ 
we have $I\ne\emptyset$.  If $|I|=1$ then the polynomial
$\prod_{\nu\in I}u_\nu$ belongs to $\Ii_0$ by assumption.
Hence assume $|I|\ge2$.   Then
\begin{equation}\label{eq:nonempty}
     \nu\in I\qquad\IMP\qquad \Delta_{\{\nu\}}\ne\emptyset.
\end{equation}
We shall prove that there exists a vector $d\in\Dd_\eff(\tau)\setminus\{0\}$ 
such that $d_\nu=1$ for $\nu\in I$ and 
$d_\nu\leq 0$ for $\nu\notin I$.  To see this,
consider the vector $\sum_{i \in I}\bar e_i\in\bar\Lambda$.
Since the union of the cones dual to vertices is the whole space 
$\bar\tt$, it follows that there exists an 
index set $J\subset\{1,\dots,n\}$ such that $|J|=n-k$,
$\Delta_J\ne\emptyset$, and
$\sum_{i \in I}\bar e_i\in\check{\Delta}_J$. 
Hence there exists nonnegative real numbers $c_j$ such that 
$$
     \sum_{i \in I}\bar e_i = \sum_{j\in J}c_j\bar e_j.
$$
By Lemma~\ref{le:lattice}, the set $\{\bar e_j\,|\,j\in J\}$ is an integer
basis of $\bar\Lambda$.  Hence the $c_j$ are actually integers 
and, after shrinking $J$, we may assume that $c_j>0$ for 
all $j\in J$. Define $d\in\Z^n$ by 
$$
     d_\nu := \left\{\begin{array}{rl}
     1,&\mbox{if }\nu\in I\setminus J,\\
     -c_\nu,&\mbox{if }\nu\in J\setminus I,\\
     1-c_\nu,&\mbox{if }\nu\in I\cap J,\\
     0,&\mbox{if }\nu\notin I\cup J.
     \end{array}\right.
$$
Then $\sum_{\nu=1}^nd_\nu \bar e_\nu=0$ and hence, by
Lemma~\ref{le:lattice}, there exists a lattice vector 
$\lambda\in\Lambda$ such that 
$d_\nu=\la\w_\nu,\lambda\ra$ for $\nu=1,\dots,n$.
By~(\ref{eq:nonempty}), we have $\Delta_{\{\nu\}}\ne\emptyset$ 
for every $\nu\in I\cup J$. Hence $d_\nu=0$ whenever 
$\Delta_{\{\nu\}}=\emptyset$, and this implies $\lambda\in\Lambda(\tau)$. 

We prove that 
$d\in\Dd_\eff(\tau)\setminus\{0\}$.
Let $\phi:\bar\tt\to\R$ be the support function of $\Delta$ as in 
Appendix~\ref{app:convex}. By property~(P2), we have 
$$
     \sum_{i\in I}\phi(\bar e_i) 
     \le \phi\left(\sum_{i\in I}\bar e_i\right)  
     = \phi\left(\sum_{j\in J}c_j\bar e_j\right)
     = \sum_{j\in J}c_j\phi(\bar e_j).
$$
Here the last equation follows from property~(F2) and the fact that
the set $\{\bar e_j\,|\,j\in J\}$ spans the cone 
$\check{\Delta}_J$. Now, by~(\ref{eq:Delta}) and the 
definition of $\phi$, we have 
$
      \phi(\bar e_\nu) \ge\zeta_\nu,
$
with equality if and only if the face $\Delta_{\{\nu\}}$ 
is nonempty. Moreover, $d_\nu=0$ whenever $\Delta_{\{\nu\}}=\emptyset$. 
This implies 
$$
       0 \ge \sum_{\nu=1}^nd_\nu\phi(\bar e_\nu) 
       = \sum_{\nu=1}^nd_\nu\zeta_\nu
       = \sum_{\nu=1}^n\inner{\w_\nu}{\lambda}\zeta_\nu 
       = -\inner{\tau}{\lambda}.
$$
If we replace $\tau$ by another vector $\tau'$ 
in the same chamber, the fan $\Sigma$ remains the same, so the above
argument yields the same  
vector $\lambda\in\Lambda$. This shows that 
$\la\tau',\lambda\ra\ge0$ for every $\tau'$ in the chamber
of $\tau$. So $\lambda\in\Lambda_\eff(\tau)$ 
and $d\in\Dd_\eff(\tau)$.  
Since $\Delta_I=\emptyset$ and $\Delta_J\ne\emptyset$,
we have $I\ne J$ and hence $d\ne 0$. 

We prove that $I\cap J=\emptyset$. Otherwise let
$\nu_0\in I\cap J$ and $I':=I\setminus\{\nu_0\}$. 
Then $d_\nu\leq 0$ for $\nu\notin I'$.  Hence the argument 
in the proof of $\Ii_0\subset\Ii$ shows that $\Delta_{I'}=\emptyset$.
But this contradicts the minimality assumption on $I$. 
Hence $I\cap J=\emptyset$ as claimed. 
It follows that the vector $d$ satisfies $d_\nu=1$ 
for $\nu\in I$ and $d_\nu\le 0$ for $\nu\notin I$.
Since $d\in\Dd_\eff(\tau)\setminus\{0\}$ we deduce that
$$
\prod_{\nu\in I}u_\nu =u^{d^+}\in\Ii_0.
$$
This shows that $\Ii\subset\Ii_0$ and hence 
$\Ii=\Ii_0$.
\end{proof}


\end{document}